\documentclass[a4paper,10pt]{article}
\usepackage[latin1]{inputenc}
\usepackage{amsthm}
\usepackage{amsmath,amssymb,amsfonts}
\usepackage{graphicx}
\usepackage{epic}
\usepackage[french,english]{babel}
\usepackage[mathscr]{euscript}
\usepackage{enumerate}
\usepackage{xspace}

\input xy
\xyoption{all}
\usepackage[dvips,all]{xy}
\input xypic
\xyoption{v2}
\usepackage[pdftex]{hyperref}
\usepackage{url}
\usepackage{multirow}
\numberwithin{equation}{section}
\newtheorem{thm}{Théorème}[section]
\newtheorem{conj}[thm]{Conjecture}

\newtheorem{prop}[thm]{Proposition}
\newtheorem{lm}[thm]{Lemme}
\newtheorem{cor}[thm]{Corollaire}

\theoremstyle{definition}
\newtheorem{defi}[thm]{D\'efinition}

\newtheorem*{cons}{Conséquence}

\newtheorem{ex}[thm]{Exemple}
\newtheorem{Rq}[thm]{Remarque}

\theoremstyle{remark}

\newtheorem*{mercis}{Remerciements}

\theoremstyle{plain}
\newtheorem{mainthm}{Th\'eor\`eme}

\newtheorem{mainconj}[mainthm]{Conjecture}

\newcommand{\Z}{\mathbb{Z}}
\newcommand{\N}{\mathbb{N}}
\newcommand{\K}{\mathbb{\Bbbk}}

\newcommand{\g}{\mathfrak{g}}
\newcommand{\h}{\mathfrak{h}}

\newcommand{\kk}{\mathfrak{k}}

\newcommand{\pp}{\mathfrak{p}}
\newcommand{\qq}{\mathfrak{q}}

\newcommand{\lf}{\mathfrak{l}}

\newcommand{\sfr}{\mathfrak{s}}
\newcommand{\wfr}{\mathfrak{w}}

\newcommand{\CC}{\mathfrak{C}}
\newcommand{\Cnil}{\mathfrak C^{\textrm{nil}}(\g)}
\newcommand{\cnil}{\mathfrak C^{\textrm{nil}}(\pp)}
\newcommand{\Od}{\mathcal{O}}
\newcommand{\NN}{\mathcal{N}}

\newcommand{\R}{\mathbb{R}}

\newcommand{\lnq} {<}
\newcommand{\gnq} {>}
\newcommand{\sld}{\mathfrak{sl}_2}
\newcommand{\sln}{\mathfrak{sl}}
\newcommand{\so}{\mathfrak{so}}
\newcommand{\spn}{\mathfrak{sp}}
\newcommand{\Striplet}{\mathcal{S}\textrm{-triplet}}
\newcommand{\ad}{\mathrm{ad }}
\newcommand{\Ad}{\mathrm{Ad }}
\newcommand{\pr}{\mathrm{pr}}
\newcommand{\rk}{\mathrm{rk }}
\newcommand{\codim}{\mathrm{codim}}

\newcommand{\gl}{\mathfrak{gl}}

\newcommand{\Phit}{\tilde{\Phi}}
\newcommand{\Id}{\mathrm{Id}}

\newcommand{\aut}{\mathfrak{aut}}
\newcommand{\GL}{\mathrm{GL}}
\newcommand{\SL}{\mathrm{SL}}
\newcommand{\Nej}{\mathcal N_e^{(j)}}
\newcommand{\Cej}{\CC(e)^{(j)}}
\newcommand{\Ceij}{\CC(e_i)^{(j)}}
\newcommand{\Lie}{\mathfrak{Lie}}
\title{Composantes irréductibles de la variété commutante nilpotente d'une algèbre de Lie symétrique semi-simple}
\author{{\sc
    Micha\"el~Bulois}\thanks{{\url{Michael.Bulois@univ-brest.fr}}}
  \\ 
  D\'epartement de math\'ematiques \\ Universit\'e de Brest \\
  29238 Brest cedex~3, France}
\begin{document}
\selectlanguage{french}
\date{}
\maketitle
\begin{abstract}
Soit $\theta$ une involution de l'algèbre de Lie semi-simple de dimension finie $\g$ et $\g=\kk\oplus\pp$ la décomposition de Cartan associée. 
La variété commutante nilpotente de l'algèbre de Lie symétrique $(\g,\theta)$ est formée des paires d'éléments nilpotents $(x,y)$ de $\pp$ tels que $[x,y]=0$.
Il est conjecturé que cette variété est équidimensionnelle et que ses composantes irréductibles sont indexées par les orbites d'éléments $\pp$-distingués. 
Cette conjecture a été démontrée par A.~Premet dans le cas $(\g\times\g,\theta)$ avec $\theta(x,y)=(y,x)$. Dans ce travail, nous la prouvons dans un grand nombre d'autres cas.
\end{abstract}

{\selectlanguage{english}
\begin{abstract}
Let $\theta$ be an involution of the finite dimensional semisimple Lie algebra $\g$ and $\g=\kk\oplus\pp$ be the associated Cartan decomposition. 
The nilpotent commuting variety of $(\g,\theta)$ consists in pairs of nilpotent elements $(x,y)$ of $\pp$ such that $[x,y]=0$. 
It is conjectured that this variety is equidimensional and that its irreducible components are indexed by the orbits of $\pp$-distinguished elements.
This conjecture was established by A.~Premet in the case $(\g\times\g,\theta)$ where $\theta(x,y)=(y,x)$.
In this work we prove the conjecture in a significant number of other cases.
\end{abstract}
}

\section*{Introduction}
Soit $\g$ une algèbre de Lie réductive, de dimension finie, définie sur un corps $\K$ algé\-bri\-que\-ment clos de caractéristique zéro. Soit
$G$ son
groupe adjoint de sorte que $\g=\Lie(G)$. Soit $\theta$ un automorphisme involutif de $\g$ et $\g=\kk\oplus\pp$ la
décomposition associée en $\theta$-espaces propres de valeurs propres respectives $+1$ et $-1$. Ceci nous donne une algèbre
de Lie symétrique $(\g,\kk)$. Notons $K$ le groupe adjoint de $\kk$. Le travail de R.W.~Richardson \cite{Ri} a montré que
$\CC(\g)=\{(x,y)\in\g\times\g\mid[x,y]=0\}$, la variété commutante de $\g$, est irréductible. Suivant une conjecture \cite{Ba} de V.~Baranovsky, A.~Premet
a montré \cite{Pr} que $\Cnil$, la variété
commutante nilpotente de $\g$ est équidimensionnelle et a indexé ses composantes irréductibles par les orbites
distinguées.
Par ailleurs, l'irréductibilité de la variété commutante de $\pp$, $\CC(\pp)$, a été étudiée dans
 \cite{PY, Pa1, Pa2, Pa3, SY1, SY2}. Nous nous intéressons à $\cnil$, la variété commutante
nilpotente de $\pp$. 
Soit $\NN$ le c\^one des éléments nilpotents de $\pp$,
on note 
$$
\cnil = \bigl\{(x,y) \in \NN \times \NN : [x,y]
= 0\bigr\} = \Cnil \cap (\pp \times \pp)
$$
la variété commutante nilpotente de $\pp$.  Le but de ce travail
est d'\'etablir, dans beaucoup de cas, la conjecture
suivante. (Rappelons qu'un élément de $\pp$ est dit
$\pp$-distingu\'e si son centralisateur dans $\pp$ ne contient
que des éléments nilpotents.)
\begin{mainconj}
\label{conjec}
La variété $\cnil$ est \'equidimensionnelle de dimension
$\dim \pp$. Ses composantes irr\'eductibles sont index\'ees
par les orbites d'éléments $\pp$-dis\-tin\-gués. 
\end{mainconj}

Il est facile de voir qu'il suffit de prouver le r\'esultat
lorsque la paire sym\'etrique $(\g,\kk)$ est
irr\'eductible. 
En
adoptant les notations de \cite[p.~518]{He}, nous
\'etablissons la conjecture dans les cas AIII, CII, DIII,
E$\{$II--IX$\}$, FI, FII et GI.  Dans les autres cas nous
obtenons des r\'esultats qui renforcent sa validit\'e.

\medskip

La nature des composantes irr\'eductibles potentielles de
$\pp$ peut \^etre comprise \`a l'aide de
consid\'erations g\'en\'erales contenues dans la
première section. Certaines des m\'ethodes utilis\'ees
g\'en\'eralisent celles de \cite{Pr}; on est en particulier
amen\'e \`a introduire la notion d'élément presque
$\pp$-distingu\'e (cf.~\ref{ppdd}). Un tel élément $e$, s'il n'est
pas $\pp$-distingu\'e, d\'efinit une variété $\CC(e)$
susceptible d'\^etre une composante, dite \'etrange, de
dimension $< \dim \pp$ dans $\cnil$. L'essentiel du
travail consiste donc \`a montrer que les variétés de la forme
$\CC(e)$ pour de tels $e$ ne fournissent pas de
composantes irr\'eductibles.

Les sections~2 et~3 donnent une classification des
éléments presque $\pp$-distingu\'es en termes
d'($ab$-)diagammes de Young (cf.~\cite{Ot2,Ot1}) dans le cas o\`u
$\g$ est classique. Cette classification permet de prouver
la conjecture dans les cas AIII, CII et DIII.

Dans les sections~4 et~5 on montre qu'un certain nombre de
composantes \'etranges ne peuvent appara\^{\i}tre dans les
cas AI, AII, CI et BDI. Ces r\'esultats assurent, par
exemple, que la conjecture est vraie dans les cas AI en rang
$\leqslant 4$, AII en rang $\leqslant 3$, CI en rang $\leqslant 7$ et BDI en
rang $\leqslant 2$.

La section~6 traite du cas o\`u $\g$ est exceptionnelle. \`A
l'aide de tables \'etablies par D.~Z.~Djokovic on y
d\'emontre la conjecture dans tous les cas, sauf celui de EI
o\`u deux composantes \'etranges restent \`a \'eliminer. 

L'appendice~7 est un complément permettant de décrire une 
classe d'éléments appelés $\pp$-self-large. Ces éléments 
sont intimement liés à la méthode issue de la section~1.4
visant à éliminer un certain nombre de composantes étranges.
Il fait suite à l'article de D.~Panyushev \cite{Pa4} qui traite 
du cas des algèbres de Lie.

\begin{mercis}
D.~Panyushev nous a indiqu\'e qu'il a obtenu des r\'esultats
semblables aux notres. Nous le remercions
de nous en avoir inform\'e.
Je tiens également à remercier le rapporteur pour ses très pertinentes
remarques et suggestions. 
\end{mercis}

\section{Généralités}
Rappelons quelques résultats tirés de \cite{KR}. Tout élément $t\in \pp$ s'écrit de façon unique $t=s+n$ où $s$ et $n$
sont des éléments de $\pp$
respectivement semi-simple et nilpotent (via l'action adjointe $\ad_{\g}$ sur $\g$). On appelle tore de dimension
$r$ toute algèbre de Lie commutative constituée d'éléments semi-simples; on notera souvent $T_r$ un tel tore. On appelle
rang de
$\pp$ et on note $\rk(\pp)$ la dimension commune des tores maximaux de $\pp$. L'ensemble $\NN$ est un cône; qui est
une variété équidimensionnelle de dimension $\dim \pp-\rk(\pp)$.
Le cône $\NN$ est stable sous l'action de $K$ et se
décompose en un nombre fini d'orbites. On notera $\Od(e)$ la $K$-orbite d'un élément nilpotent $e$.  
Si $X\subset\g$ est une partie quelconque de $\g$, on note $\NN(X)$ l'ensemble des éléments nilpotents de $\g$
contenus dans $X$. Pour $x\in \g$ on pose $X^x=\{y\in X\mid [x,y]=0\}$.
Rappelons aussi que pour tout élément $e\in \pp$, on a \cite[Proposition 5]{KR} $$\dim \Od(e)=\dim\kk-\dim \kk^e=\dim \pp-\dim \pp^e.$$
Enfin lorsqu'il n'y aura pas d'ambiguïté, $\pr_1$ désignera une application de projection sur la première variable.

\subsection{Réduction au cas simple}
\label{sec1}

On pose $\g'=[\g,\g]$; $\kk'=\g'\cap\kk$ et $\pp'=\g'\cap \pp$. On a alors $\g'=\kk'\oplus\pp'$, ce qui nous donne une paire symétrique semi-simple $(\g',\kk')$. Notons que d'après la définition de $\cnil$, on a $\cnil=\mathfrak C^{\textrm{nil}}(\pp')$.
Pour l'étude de $\cnil$, on peut donc supposer sans perte de généralité que $\g$ est semi-simple. Ce sera le cas dans tout le reste de l'article.


De plus, si $\g=\bigoplus_i\g_i$ est une décomposition de $\g$ en algèbres de Lie simples (cf. \cite[20.1.7]{TY}), l'action de $\theta$ envoie chaque $\g_i$ sur un $\g_{\theta(i)}$. On a alors deux cas possibles:
\begin{itemize}
\item[a)] $i=\theta(i)$, dans ce cas $\g_i$ est $\theta$-stable et $\g_i=\kk_i\oplus\pp_i$
\item[b)] $i\neq \theta(i)$; dans ce cas $\theta(\theta(i))=i$ et on peut supposer que $\g_i\times\g_{\theta(i)}\cong \g_i\times\g_i$ avec $\theta(x,y)=(y,x)$. On pose $\kk_i=\{(x,x)\mid x\in \g_i\}$ et $\pp_i=\{(x,-x)\mid x\in\g_i\}$.
\end{itemize}
On a alors $\NN(\pp)=\bigoplus_i\NN(\pp_i)$ et $\cnil=\bigoplus_i \mathfrak C^{\textrm{nil}}(\pp_i)$.

Par ailleurs, dans le cas b), on a $\pp_i{\cong}\g_i$ par $(x,-x)\mapsto x$. Cet isomorphisme envoie $\NN(\pp_i)$ sur $\NN(\g_i)$ et $\mathfrak C^{\textrm{nil}}(\pp_i)$ sur $\mathfrak C^{\textrm{nil}}(\g_i)$, dont Premet a décrit les composantes irréductibles \cite{Pr}.

La classification des composantes irréductibles de $\mathfrak C^{\textrm{nil}}(\pp_i)$ pour $\g_i$ simple, suffit donc pour obtenir la classification des composantes irréductibles de $\cnil$.

\subsection{Paramétrisation par les orbites}

Rappelons que $\NN$ désigne le cône des éléments nilpotents de $\pp$. Soit $e\in \NN$. Pour $e=\{0\}$ on pose $\g(e,0)=\g$ et $\g(e,i)=0$ si $i\in\Z^*$. 
Si $e\neq0$, il existe un $\Striplet$ normal $(e,h,f)$ contenant $e$ (cf. \cite[Proposition 4]{KR}).
On pose $\g(i,h)=\{x\in \g\mid [h,x]=ix\}$ et on a $\g=\bigoplus_{i\in\Z} \g(i,h)$ (cf. \cite[19.2.7]{TY}). Comme $h\in \kk$, cette décomposition est $\theta$-stable et on pose $\kk(i,h)=\g(i,h)\cap \kk$ et $\pp(i,h)=\g(i,h)\cap \pp$. On pose, par ailleurs, $\g(e,i)=\g(i,h)\cap\g^e$. 
Comme $e\in \pp$, $\g(e,i)$ est $\theta$-stable et on note $\pp(e,i)=\g(e,i)\cap\pp$. On a alors 
$$\pp^e=\bigoplus_{i\in \N}\pp(e,i).$$

On sait que $\g(e,0)$ est le stablisateur de $(e,h,f)$ dans $\g$, c'est donc une sous-algèbre réductive
dans $\g$ (cf. \cite[20.5.13]{TY}). L'ensemble de ses éléments nilpotents est donc $\NN(\g)\cap\g(e,0)$ et on a
$$\NN(\g)\cap\g^e=\NN(\g(e,0))\times \bigoplus_{i\gnq 0}\g(e,i).$$
Par la discussion ci-dessus sur la $\theta$-stabilité, on a de plus que:
$$\NN\cap\pp^e=\NN(\pp(e,0))\times \bigoplus_{i\gnq 0} \pp(e,i).$$
C'est une variété qui contient le même nombre de composantes irréductibles que $\NN(\pp(e,0))$ que l'on indexe par un ensemble $I_e$. On a alors:
$$\NN\cap \pp^e=\bigcup_{j\in I_e} \Nej \qquad\textrm{où}\qquad \Nej=(\NN(\pp(e,0))_j\times \bigoplus_{i\gnq 0} \pp(e,i).$$

\begin{defi}
On pose $$\Cej:=\overline{\Ad K.(e,\Nej)}\subset \cnil \quad \mathrm{ et }\quad\CC(e)=\overline{\Ad K.(e,\NN\cap\pp^e)}=\bigcup_{j\in I_e} \Cej.$$ On dit que $e$ engendre $\Cej$.
\end{defi}

Les sous-variétés du type $\Cej$ sont irréductibles. Or par \cite[Théorème~2]{KR}, il existe un nombre fini de $K$-orbites nilpotentes dans $\pp$, dont on notera $e_1,\dots, e_k$ des représentants, de sorte que:
$$\cnil=\bigcup_{\substack{i=1,\dots, k\\j\in I_{e_i}}} \Ceij.$$ 
Cette union étant finie, on peut en déduire que les composantes irréductibles de $\cnil$ sont de la forme $\Ceij$ pour $i\in[\![1;k]\!]$ et $j\in I_{e_i}$.

\subsection{\'Eléments $\pp$-distingués et conjecture}
La notion d'élément nilpotent $\pp$-distingué va s'avérer très importante pour la suite.
\begin{defi} Soit $e\in \NN.$ L'élément $e$ est dit $\pp$-distingué si $\pp^e\subset\NN$.
\end{defi}
\begin{lm}
L'élément $e$ est $\pp$-distingué si et seulement si $\pp(e,0)=\{0\}$, ce qui revient à dire que $\pp^e=\bigoplus_{i\gnq0} \pp(e,i)$.
\end{lm}
\begin{proof}
Si $\pp(e,0)=\{0\}$ alors $\pp^e=\bigoplus_{i\gnq0} \pp(e,i)\subset\NN$. Réciproquement, si $\pp(e,0)\neq\{0\}$ alors $\pp(e,0)$ est le $-1$-espace propre de $\g(e,0)$ qui est réductif dans $\g$. Il contient donc des éléments semi-simples non-triviaux et $e$ n'est pas distingué.
\end{proof}
Soit $e\in \NN.$ On veut calculer la dimension des $\Cej$.
On s'intéresse pour cela à l'application dominante entre variétés irréductibles 
$$\xi: \left\{\begin{array}{r c l} K \times\Nej&\rightarrow& \Cej\\
                           (g,x) &\mapsto& (g.e,g.x)\end{array}\right..$$
Pour tout $x\in \Nej$, la fibre $\xi^{-1}(\xi(1,x))$ est l'ensemble des $(g,g^{-1}.x)$ avec $g\in K^e$.
On en déduit que $\dim \xi^{-1}(\xi(1,x))= \dim K^e=\dim \kk^e$. Ceci reste vrai pour toute fibre non-vide, et en notant $p=\dim\pp$ on a (cf. \cite[Théorème 15.5.3]{TY}) 
\begin{eqnarray*}\dim \Cej&=&\dim K+\dim \Nej-\dim\kk^e\\
                            &=&p-(\dim\pp^e-\dim \Nej)\\
&=&p-\codim_{\pp(e,0)} \;\,\Nej\\
&=&p-\rk(\pp(e,0)) \end{eqnarray*}
où la dernière égalité découle de \cite[Théorème 3]{KR}. On remarque que les sous-variétés $\CC(e)$ sont équidimensionnelles.

\begin{defi}
Pour un élément nilpotent $e\in \pp$, on définit le défaut de $e$ comme étant le rang de $\pp(e,0)$. On le note $\delta(e)$. Le défaut est invariant sous l'action de $K$. 
\end{defi}

\begin{prop}\label{dimppreconj}
La variété $\cnil$ est de dimension $p$ et ses composantes irréductibles de dimension maximale sont les $\CC(e)$ avec $e$ représentant d'orbite $\pp$-distinguée. 
\end{prop}
\begin{proof}
Par le calcul précédent, on voit que $\dim \Cej=p$ si et seulement si $e$ est de défaut nul, ce qui est équivalent à dire que $e$ est $\pp$-distingué. 
Dans le cas contraire, on a $\dim \Cej\lnq p$.
Par ailleurs, comme les éléments nilpotents $\pp$-réguliers sont $\pp$-distingués, il existe des éléments $\pp$-distingués ce qui suffit à montrer que $\dim\cnil=p$. 
Maintenant, si $e$ est $\pp$-distingué, $\NN\cap\pp^e=\pp^e$ donc $\CC(e)$ est irréductible.
On a de plus $\pr_1(\CC(e))=\overline{K.e}$. On en déduit que les $\CC(e_{i})$ sont distincts pour des éléments $e_{i}$ appartenant à des $K$-orbites différentes.
\end{proof}

On appellera \emph{composante étrange} une composante irréductible de $\cnil$ de dimension strictement inférieure à $p$. 

\begin{Rq}\label{rqconj} En rang 0, \emph{i.e.} si $\pp=\{0\}$ ou encore si $\theta$ est triviale, on a $\cnil=\{0\}\times\{0\}=\CC(0)$.
\item
En rang 1, le commutant d'un élément $x\in \pp$ est $\K x$. Tous les éléments nilpotents sont donc $\pp$-distingués et les composantes irréductibles de $\cnil$ sont les $\CC(e)$ où $e$ est nilpotent. 
\end{Rq}

\subsection{\'Eléments presque $\pp$-distingués}
Nous allons maintenant identifier les composantes étranges en introduisant la notion d'élément presque $\pp$-distingué, 
qui généralise celle d'élément presque distingué (cf. \cite{Pr}). Les résultats qui suivent, et particulièrement la proposition \ref{pd}, 
sont largement inspirés de \cite[Proposition 2.1]{Pr}.

\begin{defi}\label{ppdd}
On dit que $e\in \NN$ est presque $\pp$-distingué si $\pp(e,0)$ ne contient pas d'élément nilpotent non nul, c'est à dire si $\pp(e,0)$ est un tore $T_{\delta(e)}$. C'est équivalent au fait que $\delta(e)=\dim \pp(e,0)$.
\end{defi}

\begin{lm}\label{orbite}
Si $e\in \NN$ et $j\in I_e$ sont tels que $\Cej$ est une composante irréductible de $\cnil$, alors $\Nej\subseteq\overline{\Od(e)}$.
\end{lm}
\begin{proof} 
Remarquons tout d'abord que le groupe $\GL(2)$ agit sur $\pp\times\pp$ via:
$$\left(\begin{array}{c c} \alpha & \beta \\ \gamma & \delta\\ \end{array}\right) . (x,y)=(\alpha x +\beta y, \gamma x+\delta y).$$
Comme toute combinaison linéaire d'éléments commutant dans $\NN$ est encore dans $\NN$, la variété $\cnil$ est $\GL(2)$-invariante. Puisque $\GL(2)$ est connexe, il fixe la composante irréductible $\Cej$. En particulier $\Cej$ est stable sous l'application $\sigma: (x,y)\rightarrow (y,x)$ de $\pp\times\pp$. 
Comme $\pr_1(K.(e,\Nej))=\Od(e)$,  on a $\pr_1(\Cej)\subseteq\overline{\Od(e)}$. Ceci implique en particulier que
$$\Nej= (\pr_1\circ\sigma)(e,\Nej)\subseteq \overline{\Od(e)}.$$
\end{proof}

Notons ici une remarque liée à la preuve précédente: si $\Nej\subseteq\overline{\Od(e)}$ alors $\GL(2).\Cej=\Cej$. On ne peut donc pas obtenir plus de résultats en considérant l'action de $\GL(2)$ à la place de celle de $\sigma$. 

\begin{lm}\label{norbite}
On a $\bigoplus_{i\geqslant2} \pp(e,i) \subset \overline{\Od(e)}$. 
\end{lm}
\begin{proof}
On peut supposer $e\neq0$. On considère la sous-algèbre parabolique $\qq=\bigoplus_{i\geqslant 0} \kk(i,h)$ de l'algèbre de Lie $\kk$. D'après (\cite[29.4.3]{TY}), il existe un sous-groupe $Q$ de $K$ ayant $\qq$ pour algèbre de Lie. On a alors $[\qq,e]=\bigoplus_{i\geqslant 2} \pp(i,h)$, donc $Q.e$ est une sous-variété de $\bigoplus_{i\geqslant 2} \pp(i,h)$ de même dimension. On en déduit que 
$\overline{Q.e}=[\qq,e]=\bigoplus_{i\geqslant 2} \pp(i,h)$. L'affirmation du lemme est alors immédiate.
\end{proof}

Pour éliminer des possibilités de composantes étranges, nous allons utiliser le lemme \ref{orbite}. 
D'après le lemme \ref{norbite}, cela ne peut éventuellement s'appliquer que dans le cas où $\pp(e,0)$ ou $\pp(e,1)$ contient des éléments nilpotents non nuls.

\begin{prop}\label{pd}
Les composantes irréductibles de $\cnil$ sont de la forme $\CC(e)$ avec $e$ presque $\pp$-distingué. 
\end{prop}
\begin{proof}
Choisissons $e$ et $j\in I_e$ tels que $\Cej$ soit une composante irréductible de $\Cnil$ et supposons que $e$ ne soit pas presque $\pp$-distingué. 
Alors $\NN\cap\pp(e,0)$ est non-trivial et $\NN(\pp(e,0))_j$ contient un élément non-nul $e_0$. 
Par définition de $\Nej$, on a $e+e_0\in \Nej$, donc $e+e_0\in \overline{\Od(e)}$ par le lemme \ref{orbite}.
Si $e=0$ on a une contradiction, on suppose donc $e\neq0$.

Comme $\g(e,0)$ est réductive dans $\g$, on peut inclure $e_0$ dans un $\Striplet$ normal $(e_0,h_0,f_0)$ de $\g(e,0)$. 
Soit $\sfr_0$ l'algèbre de Lie de dimension trois engendrée par $(e_0,h_0,f_0)$. Comme $\sfr_0\subseteq \g(e,0)$, $(e+e_0,h+h_0,f+f_0)$ est un $\Striplet$ normal.
Maintenant, $e\in \g(h_0,0)$, $e_0\in \g(h_0,2)$ donc l'application $\tau_\lambda\in \K$ ($\lambda\in \K^*$) définie dans \cite[38.6.2]{TY} envoie $e+e_0$ sur $e+\lambda^2e_0$. 
On en déduit que $e+\K^*e_0\subset \Od(e+e_0)$ et donc $e\in \overline{\Od(e+e_0)}$. Rappelons que d'après le paragraphe 
précédent, on a $e+e_0\in\overline{\Od(e)}$. On en déduit que $e+e_0$ et $e$ sont $K$-conjugués. 
Ceci implique que $h$ et $h+h_0$ sont également $K$-conjugués. On va maintenant chercher une contradiction.

On a $L(h,h_0)=L(h,[e_{0},f_{0}])= L([e_{0},h],f_{0})=0$ où $L$ désigne la forme de Killing de $\g$. 
Ceci implique que $L(h+h_0,h+h_0)=L(h,h)+L(h_0,h_0)$. Mais $h$ et $h+h_0$ sont conjugués et $L$ est invariante sous l'action de $G$, d'où $L(h_0,h_0)=0$. 
Or, puisque $h_0$ est l'élément semi-simple d'un $\Striplet$, l'action adjointe de $h_{0}$ est à valeurs propres entières, et donc $L(h_0,h_0)=0$ si et seulement si $h_0=0$. 
Ceci est impossible par le choix de $e_0$. On a donc montré par l'absurde que $e$ est presque $\pp$-distingué.

Maintenant, comme $\pp(e,0)$ ne contient pas d'éléments nilpotents et que $\NN\cap\pp^e=\bigoplus_{i\gnq0} \pp(e,i)$ est irréductible, on en déduit que $\CC(e)=\Cej$ est la composante irréductible de $\cnil$ engendrée par $e$.
\end{proof}
\begin{Rq}Notons que d'après la preuve de la proposition précédente $\CC(e)$ est irréductible pour $e$ presque $\pp$-distingué.
On a donc $\CC(e)=\overline{K.(e,\NN\cap\pp^e)}$ lorsque $e$ est presque $\pp$-distingué et la présence des indices $j\in I_{e}$ n'est plus nécessaire. 
\end{Rq}
\begin{cor}\label{corpdd}
Les composantes étranges sont engendrées par des éléments presque $\pp$-distingués non $\pp$-distingués. 
\end{cor}
D'après la proposition \ref{dimppreconj} et le corollaire \ref{corpdd} la conjecture \ref{conjec} 
est équivalente au fait que $\cnil$ ne comporte pas de composantes étranges,  ou encore:
\begin{conj}\label{conj}
Si $e_{1}$ est presque $\pp$-distingué, il existe un élément $e_{2}$ $\pp$-distingué tel que $\CC(e_{1})\subseteq\CC(e_{2})$.
\end{conj}
D'après les résultats de la section \ref{sec1}, cette conjecture est vraie si et seulement si elle est vraie pour toute algèbre de Lie symétrique simple.

\section{Centralisateurs}\label{centralisateurs}
Le but de cette section est de caractériser les centralisateurs de $\Striplet$s dans le cas classique pour obtenir une classification des éléments presque $\pp$-distingués.
En particulier nous allons retrouver que les éléments de défaut nul (\emph{i.e.} $\pp$-distingués) correspondent bien au éléments compacts décrits dans \cite{PT}. 
Le travail de \cite{PT} donne une description des éléments $\pp$-distingués et utilise le fait que les orbites d'éléments compacts dans le cas réel correspondent, 
via la correspondance de Kostant-Sekiguchi, aux orbites d'éléments $\pp$-distingués. Notre démarche pour décrire les centralisateurs lui sera en fait assez similaire.
Dans toute la section, $V$ désignera un espace vectoriel de dimension $n$. On supposera que $\g\subset\sln(V)$ est une algèbre de Lie simple et que, contrairement à l'introduction, $G\subset \SL(V)$
est le plus petit groupe algébrique dont l'algèbre de Lie contient $\g$ (cf. \cite[Chapitre~24]{TY}). Cette modification n'a pas d'incidence sur les orbites de $\g$.
On se restreint aux cas classiques, c'est à dire que l'on est dans l'une des situations suivantes:
\begin{itemize}
\item Type A: $\g=\sln(V)$ auquel cas $G=SL(V)$.
\item Type BD: $\g=\so_{\Phi}(V)=\aut(V,\Phi)$ est l'algèbre de Lie orthogonale stabilisant une forme bilinéraire symétrique non dégénérée $\Phi$. Dans ce cas $G=SO_{\Phi}(V)$.
\item Type C: $\g=\spn_{\Phi}(V)=\aut(V,\Phi)$ est l'algèbre de Lie symplectique stabilisant une forme bilinéraire alternée non dégénérée $\Phi$. Dans ce cas $G=SP_{\Phi}(V)$.
\end{itemize}

\begin{defi} \label{cas}\cite[Théorème 3.4]{GW} Enumérons l'ensemble des involutions de Cartan $\theta$ définies sur $\g$ dans chacun de ces types. 
\begin{itemize}
\item Type $A$: 
\begin{enumerate}
\item[(AI)] $\theta(x)=-T {}^tx {}^tT$ pour $x\in \g$ et où $T=T^{-1}={}^tT$. Ceci est équivalent au fait que $\kk=\aut(V,\Phi)$ où $\Phi(x,y)={}^txTy$ est la forme bilinéaire symétrique associée à $T$. La propriété $T=T^{-1}={}^tT$ définit de façon unique $\theta$ (à conjugaison près dans $G$).
\item[(AII)] $\theta(x)=-T {}^tx {}^tT$ pour $x\in \g$ et où $T=T^{-1}=-{}^tT$. Ceci est équivalent au fait que $\kk=\aut(V,\Phi)$ où $\Phi={}^txTy$ est la forme bilinéaire antisymétrique associée à $T$. La propriété $T=T^{-1}=-{}^tT$ définit de façon unique $\theta$ (à conjugaison près dans $G$).
\item[(AIII)] $\theta(x)=JxJ^{-1}$ pour $x\in \g$ et où $J^2=\Id_V$. On définit
$$V_{\pm}=\{v\in V \mid Jv=\pm v\}.$$
Alors $V=V_{+}\oplus V_{-}$ et le nombre $\dim(V_{+})$ définit $\theta$ de façon unique (à conjugaison près dans $G$).
\end{enumerate}
\item Types $BD$ et $C$:
\begin{enumerate}
\item[(BDI),(CII)] $\theta(x)=JxJ^{-1}$ pour $x\in \g$ où $J$ préserve la forme $\Phi$ et vérifie $J^2=~\Id_V$. On définit
$$V_{\pm}=\{v\in V \mid Jv=\pm v\}.$$
Alors $V=V_{+}\oplus V_{-}$, la restriction de $\Phi$ à $V_{\pm}$ est non-dégénérée et le nombre $\dim(V_{+})$ définit $\theta$ de façon unique (à conjugaison près dans $G$).
\item[(DIII),(CI)] $\theta(x)=JxJ^{-1}$ pour $x\in \g$ où $J$ préserve la forme $\Phi$ et vérifie $J^2=~-\Id_V$. On définit
$$V_{\pm i}=\{v\in V \mid Jv=\pm i v\}.$$
Alors $V=V_{+}\oplus V_{-}$, la restriction de $\Phi$ à $V_{\pm i}$ est nulle et $V_{+i}$ dual à $V_{-i}$ par rapport à $\omega$. De plus l'involution $\theta$ est uniquement déterminée (à conjugaison près dans $G$).
\end{enumerate}
\end{itemize}
\end{defi}

Etant donnée une telle algèbre de Lie symétrique $(\g,\kk)$ ou $(\g,\theta)$ et $e\in\pp$ nilpotent non nul, on fixe un $\Striplet$ $(e,h,f)$ normal qui engendre une sous-algèbre de Lie de dimension 3 que l'on note $\sfr$. 
Les centralisateurs $\g^{\sfr}$ sont connus, cf. \cite{SS} ou \cite{Ja}. 
Afin d'obtenir des informations précises sur les paires symétriques $(\g^{\sfr},\theta_{\mid\g^{\sfr}})$, nous allons rappeler la description des sous-algèbres $\g^\sfr$. 

Pour $\lambda\in \Z$ notons:
$$V(\lambda)=\{v\in V\mid h.v=\lambda v\}$$
de sorte que $V=\bigoplus_{\lambda\in\Z}V(\lambda)$. 
Il existe (à isomorphisme près) un unique $\sld$-module irréductible de dimension $d$. Notons-le $\rho_d$.
Alors $V=\bigoplus_{d\in\N^*} V_d$ où $V_d$ est isomorphe à $(\rho_1)^{m_d}\otimes_{\sfr} \rho_d$ avec $\otimes_{\sfr}$ désignant le produit tensoriel de $\sfr$-modules. Ceci nous donne une décomposition 
\begin{equation}\label{dec}V_d=\bigoplus_{j\in [\![0,d-1]\!]} V(d-1-2j)\cap V_d.\end{equation}
On définit alors $V_{j,d}:=V(d-1-2j)\cap V_d$ et $m_d:=\dim V_{j,d}$, ce dernier étant un entier indépendant de $j$.
On en déduit une partition de $n$ donnée par $(d^{m_d})_{d\in \N^*}$. 
Cette partition correspond au diagramme de Young usuellement associé à la classe de conjuguaison d'un élément nilpotent de $\sln(V)$ (cf. \cite{Ot2}, par exemple).
Notons $H_d$ (resp. $H_d'$) le sous-espace vectoriel engendré par les vecteurs de plus haut (resp. bas) poids de $V_d$. Autrement dit $H_d=V_{0,d}$ et $H_d'=V_{d-1,d}$.
Soit $g$ un élément de $\gl(V)^\sfr$. Par définition il stabilise les sous espaces propres de $h$ et de $e$, et en particulier $g$ stabilise les sous-espaces 
$H_d=V(d-1)\cap \ker e$. On peut donc définir des applications restrictions indexées par $d\in \N^*$
$$\varphi_d: \left\{\begin{array}{rcl} \gl(V)^\sfr\ &\longrightarrow &\gl (H_d) \\ g &\longmapsto& g_{\mid H_d}. \end{array}\right.$$
Et on définit \begin{equation} \varphi=\bigoplus_{d\in \N^*} \varphi_d \end{equation}
On peut alors énoncer un premier résultat sur les centralisateurs:
\begin{lm}\label{cent}
Le centralisateur du $\Striplet$ $(e,h,f)$ dans $\gl(V)$ vérifie:
$$\gl(V)^\sfr\cong^\varphi\bigoplus_{d\in \N^*} \gl (H_d).$$
\end{lm}
\begin{proof}
On a vu que l'application $\varphi$ est bien définie. Nous allons dé\-mon\-trer qu'elle possède une réciproque.
Soit $(g_d)_{d\in\N^*}$ une famille d'éléments de $\bigoplus_{d\in \N^*} \gl(H_d)$, nous allons montrer qu'il existe un unique élément $g\in \gl(V)$ tel que $g$ centralise $\sfr$ et $g_{\mid H_d}=g_d$.
La formule (\ref{dec}) associée au fait que $f^j$ induit une bijection entre $H_d$ et $V_{j,d}$ nous montre qu'un tel élément $g$ est nécessairement défini uniquement sur chaque composante de $V_d$ par \begin{equation} \label{def} g_{\mid V_{j,d}}= f^j.g_d.f^{-j}. \end{equation} Comme $V$ est somme directe des $V_{j,d}$, l'élément $g$ est défini uniquement sur $V$ tout entier. Réciproquement, il est facile de vérifier qu'un élément $g$, défini par (\ref{def}) à partir d'une famille $(g_d)_{d\in\N^*}$ quelconque, commute avec $\sfr$.
\end{proof}


A partir de maintenant, on se donne une forme bilinéaire non-dégénérée $\Phi$ sur $V$ et deux éléments $\varepsilon, \eta \in\{\pm1\}$ de sorte que 
\begin{itemize}
\item $\Phi$ est symétrique ou antisymétrique ce qui se traduit par $\Phi(u,v)=\varepsilon \Phi(v,u)$ pour tout $u,v\in V$.
\item $h$ préserve $\Phi$ et $e,f$ preservent ou anti-préservent $\Phi$: c'est à dire 
$$\!\!\!\!\!\!\!\!\!\!\!\!\Phi(h.u,v)=-\Phi(u,h.v) \; ; \;  \Phi(e.u,v)=-\eta\Phi(u,e.v)\; ; \;  \Phi(f.u,v)=-\eta\Phi(u,f.v).$$
\end{itemize}
Le cas $\eta=1$ permettra de décrire $\g^{\sfr}$ quand $\g=\aut(V,\Phi)$, 
tandis que le cas $\eta=-1$ permettra de décrire $\kk^{\sfr}$ quand $\kk=\g\cap\aut(V,\Phi)$.

\begin{lm}\label{poids}
\item[(a)]Si $\lambda\neq -\mu$ alors $\Phi(V(\lambda), V(\mu))=0$. Les sous-espaces $V(\lambda)$ et $V(-\lambda)$ sont en dualité par $\Phi$.
\item[(b)]$H_d$ est dual à $H_d'=f^{d-1}(H_d)$.
\end{lm}
\begin{proof}
Soient $u\in V(\lambda), v\in V(\mu)$. 
\begin{eqnarray*}
\lambda\,\Phi(u,v)&=&\Phi(h.u,v)\\
&=& -\Phi(u,h.v)=-\mu \,\Phi(u,v),
\end{eqnarray*} 
ce qui démontre la partie (a).\\
Soit maintenant $u\in H_d$. Par (a), il existe $w\in V(-d+1)$ tel que $\Phi(u,w)\neq 0.$
Il existe $w_{1},w_{2}$ tels que $w=e.w_{1}+w_{2}$ et $w_{2}\in H_{d}'$. 
On a alors:
$$0\neq\Phi(u,w)=\Phi(u,e.w_{1}+w_{2})=-\eta \Phi(e.u,w_{1})+\Phi(u,w_{2})=\Phi(u,w_{2}).$$
Ceci prouve que la restriction de $\Phi$ à $H_{d}\times H_{d}'$ est non dégénérée.
\end{proof}

\begin{defi}
On définit $\Phit_d$ sur $H_d$ par
$$\Phit_d(u,v)=\Phi(u,f^{d-1} .v).$$
\end{defi}

\begin{lm}\label{sym}
La forme $\Phit_d$ est bilinéaire non dégénérée sur $H_d$ et elle vérifie $\Phit_d(u,v)=(-\eta)^{d-1}\varepsilon\,\Phit_d(v,u)$.
\end{lm}
\begin{proof}
L'application $f^{d-1}$ est une bijection de $H_d'$ sur $H_d$, donc par le lemme \ref{poids} (b), $\Phit_d$ est non-dégénérée sur $H_d$. La relation de symétrie est une vérification facile.
\end{proof}

Nous pouvons maintenant énoncer un second résultat sur les centralisateurs.

\begin{prop} \label{AI}
Si l'algèbre de Lie $\wfr\subseteq\sln(V)$ est l'algèbre préservant la forme $\Phi$ sur $V$ alors le centralisateur de $(e,h,f)$ dans $\wfr$ vérifie
$$\wfr^{\sfr}\cong^{\varphi}\bigoplus_{d\in\N^*} \aut(H_d,\Phit_d).$$
\end{prop}
\begin{proof}
Par le lemme \ref{cent}, on sait que tout élément de $\wfr^{\sfr}$ peut être vu comme un élément de $\bigoplus_{d\in\N^*} \gl(H_d)$.
Maintenant, le fait que $g\in \wfr$ implique que pour tout $d\in \N^*$, $(u,w)\in H_d\times H_d'$, on a $\Phi(g.u,w)=-\Phi(u,g.w)$. 
C'est à dire que pour tout $(u,v)\in H_d$, on a $0=\Phi(g.u, f^{d-1}.v)+\Phi(u,gf^{d-1}.v)=\Phit_d(g.u,v)+\Phit_d(u,g.v)$, ce qu'on peut reformuler en: $g_{\mid H_d}$ préserve $\Phit_d$.

Réciproquement, si $g$ est défini à partir d'éléments $g_d\in \aut(H_d,\Phit_d)$ et de la formule (\ref{def}), 
alors $g$ stabilise les sous-espaces $V_{j,d}$ de (\ref{dec}). 
Il suffit donc de regarder $g$ sur tous les éléments $(u,v)\in V_{j,d}\times V_{d-1-j,d}$ pour montrer que $g$ preserve $\Phi$:
\begin{eqnarray*}
\Phi(g.u,v)=\Phi(f^jg_df^{-j}.u,v)&=&(-\eta)^j\Phi(g_d f^{-j}.u, f^j.v)\\
&=& (-\eta)^j\,\Phit(g_df^{-j}.u, f^{j-(d-1)}.v)\\
&=& -(-\eta)^j\,\Phit(f^{-j}.u,g_d f^{j-(d-1)}.v)\\
&=& -(-\eta)^j\,\Phi (f^{-j}.u, f^{d-1} g_d f^{j-(d-1)}. v)\\
&=& - \Phi(u, f^{d-1-j}g_d f^{j-(d-1)}.v)=-\Phi(u,g.v).
\end{eqnarray*}
\end{proof}

Dorénavant $J$ désigne un élément de $GL(V)$ tel qu'il existe $\xi\in \{\pm 1\}$ avec $J^2=\xi\,\Id$.
Pour $g\in \gl(V)$, on note $\theta(g)=JgJ^{-1}$.
On suppose de plus que $h$ commute avec $J$, et $e,f$ anti-commutent avec $J$. C'est à dire:
$$\theta (h)=h \;;\; \theta(e)=-e \;;\; \theta(f)=-f$$
On note $\sqrt{-1}$ une racine carrée de $-1$ dans $\K$.

\begin{lm}\label{J}
\item[(a)] Le sous-espace $H_d$ est $J$-stable.
\item[(b)] Si $J$ préserve $\Phi$, alors la restriction de $(\sqrt{-1})^{d-1}J$ à $H_d$ préserve la forme~$\Phit$.
\end{lm}
\begin{proof}
Les endomorphismes $J$ et $h$ commutent, donc le $h$-sous-espace propre $V(d-1)$ est stable par $J$. De la même façon, $J$ et $e$ anti-commutent donc le noyau de $e$ est stable par $J$. On en déduit la partie (a) en remarquant que $H_d=V(d-1)\cap \ker e$.

Si $J$ préserve $\Phi$, alors pour tout $u,v\in H_d$, 
\begin{eqnarray*}\Phit((\sqrt{-1})^{d-1}J.u,(\sqrt{-1})^{d-1}J.v)&=&(-1)^{d-1}\Phi(J.u,f^{d-1}J.v)\\&=&(-1)^{2(d-1)}\Phi(J.u, Jf^{d-1}.v)\\&=&\Phi(u,f^{d-1}v)=\Phit(u,v).
\end{eqnarray*}
\end{proof}

\begin{Rq}
Pour tout $d\in \N^*$, l'involution $\theta$ est inchangée par la modification de $J$ en $(\sqrt{-1})^{d-1}J$. En effet, pour tout $g\in \gl(V)$, \\ $(\sqrt{-1})^{d-1}J g ((\sqrt{-1})^{d-1}J)^{-1}=JgJ^{-1}$. Par contre, ce qui change ce sont les formes preservées par $J$ et $(\sqrt{-1})^{d-1}J$. 
\end{Rq}

\begin{cor}\label{AIII}
Si $\gl^{\theta}(V)\subseteq \gl(V)$ est l'algèbre de Lie reductive correspondant au $1$-espace propre de l'automorphisme $\theta$ alors 
$$(\gl^{\theta})^{\sfr}\cong^{\varphi}\bigoplus_{d\in \N^*} (\gl^{\theta_d})(H_d)$$
où $\theta_d(h)$ pour $h\in H_d$ est donné par $J_{\mid H_d}hJ^{-1}_{\mid H_d}$.
\end{cor}

\begin{proof}
Le seul point restant non-trivial consiste à vérifier qu'un élément $g$ déterminé à partir d'éléments $g_d$ comme dans la formule (\ref{def}) est bien invariant par $\theta$. En effet, si $u\in V_{j,d}$:
$$\theta(g).u=\theta (f^jg_df^{-j}).u=(-1)^{j-j} f^j \theta(g_d)f^{-j} .u=g.u$$
\end{proof}

Enfin, cette dernière proposition nous permettra de traiter les derniers centralisateurs classiques envisageables. 

\begin{prop}\label{BDI}
Si $J$ preserve $\Phi$ et que $\kk=\aut(V,\Phi)^J\subseteq\gl(V)$ est l'algèbre de Lie réductive donnée par intersection de $\gl^{\theta}(V)$ et de $\aut(V,\Phi)$, alors
$$\kk^s\cong^{\varphi}\bigoplus_{d\in \N^*}(\aut(H_d, \Phit))^{(\sqrt{-1})^{d-1}J}.$$
\end{prop}
\begin{proof}
Il suffit de combiner le corollaire \ref{AIII} et la proposition \ref{AI}.
\end{proof}

\section{Classification des éléments presque $\pp$-distingués dans le cas classique}
Reprenons les notations de la section \ref{centralisateurs} où $\g\subseteq\sln(V)=\sln_n$ désignera une algèbre de Lie classique simple (cf. \ref{cas}), $\theta$ une involution non-triviale sur $\g$ et $\kk, \pp$ les sous-espaces propres de $\theta$ associés respectivement aux valeurs propres $+1$ et $-1$.
On supposera toujours que si $e\in \pp$ est un élément nilpotent non nul, on l'inclut dans un $\Striplet$ normal $(e,h,f)$. 
On laisse au lecteur le soin de traduire les assertions quand $e=0$. 
Rappelons, cf. section précédente, 
qu'à $e$ est associée canoniquement la partition $(d^{m_d})_{d\in \N^*}$ de $n$. 
On notera ces partitions sous la forme d'un diagramme de Young, qu'on appellera diagramme de Young associé à $e$.
Dans le cas où une matrice $J$ intervient dans la définition de l'involution sur l'algèbre de Lie (cf. \ref{cas}), 
on introduit la notion d'$ab$-diagramme de Young associé à $e$ (cf. par exemple \cite{Ot2}) de la façon suivante. 
On décompose $V$ en $J$-sous-espaces propres $V_a\oplus V_b$ et on fixe une base diagonalisant $J$ 
$$\{\alpha_{d,i}^j\mid d\in \N^*, i\in [\![0,d-1]\!], j\in [\![1,m_d]\!]\}$$
 adaptée au $\Striplet$, c'est à dire telle que $(\alpha^j_{d,0})_j$ forme une base de $H_d'$ et $e(\alpha^j_{d,i})=\alpha^j_{d,i+1}$. On remplit alors les cases du diagramme de Young associé à $e$ par $a$ ou $b$ suivant que $\alpha_{d,i}^j\in V_a$ ou $V_b$. Comme $\theta(e)=-e$, sur une ligne du diagramme de Young, on a alternance de $a$ et de $b$. 
On définit aussi $H_{d,a}'=V_a\cap H_d'$ (resp. $H_{d,b}'=V_b\cap H_d'$) et $a_d=\dim H_{d,a}'$ (resp. $b_d=\dim H_{d,b}'$) de sorte que $a_d+b_d=m_d$ 
et que $a_d$ (resp $b_d$) désigne le nombre de lignes de longueur $d$ commençant par un $a$ (resp. $b$).
Les éléments $(a_d)_{d\in \N^*}$ et $(b_d)_{d\in \N^*}$ sont des invariants de la $K$-classe de conjugaison de $e$.
Lorsque $J$ n'intervient pas n'intervient pas dans la définition de l'involution (type AI et AII), 
le diagramme de Young ou la donnée de $(m_{d})_{d\in\N^*}$ sera un invariant suffisant.
\begin{defi}
On appelle ($ab$-)diagramme de Young, soit un diagramme de Young (types AI et AII), soit un $ab$-diagramme de Young (autres types classiques) associé à une orbite 
$\Od(e)$ comme ci-dessus. 
\end{defi}
Rappelons qu'en vertu de la section précédente, on a des décompositions $\theta$-stables, indexées par $d$ grâce à l'application $\varphi$, des algèbres de Lie $\gl(V),\so(V)$ et $\spn(V)$. La définition suivante utilise ces décompositions.
\begin{defi}\label{delta}
On appelle défaut de la longueur $d$ et on note $\delta^e(d)$ le rang de la paire symétrique $(\g_d',\kk_d):=(\varphi_d(\g'),\varphi_d(\kk))$ où $\g'=\g$ dans les cas DIII, BDI, CI et CII et $\g'=\gl(V)$ dans les cas AI, AII et AIII. 
\end{defi}
\begin{Rq} \label{delta0}
Sous les notations précédentes, on a $\delta(e)=\sum_d \delta^e(d)$ dans les cas DIII, BDI, CI et CII et $\delta(e)=\sum_d \delta^e(d)-1$ dans les cas AI, AII et AIII. 
\end{Rq}

Les sous-sections 3.1, 3.2 et 3.3 donnent les ($ab$-)diagrammes de Young des éléments $\pp$-distingués et presque $\pp$-distingués
pour chacun des cas présentés dans la définition 2.1.
On rappelle que dans chaque cas simple classique, les \mbox{($ab$-)}diagrammes possibles sont les concaténations des ($ab$-)diagrammes primitifs que nous récapitulons 
dans le tableau suivant tiré des travaux de T.~Otha. Rappelons aussi que les $(ab$-)diagrammes sont classiquement ordonnés de façon décroissante de la plus grande ligne en haut vers la plus petite en bas et qu'on les regarde  à permutation près de lignes de même longueur.
\begin{center}
\begin{tabular}{|c| c|}
\hline
Type & ($ab$-)diagrammes primitifs\\
\hline
AI & \begin{picture} (60,15)
\put(0,0){\line(0,1){10}}
\put(0,0){\line(1,0){20}}
\put(23,0){\line(1,0){4}}
\put(33,0){\line(1,0){4}}
\put(40,0){\line(1,0){10}}
\put(50,0){\line(0,1){10}}
\put(0,10){\line(1,0){20}}
\put(23,10){\line(1,0){4}}
\put(33,10){\line(1,0){4}}
\put(40,10){\line(1,0){10}}
\put(40,3){\line(0,1){4}}
\put(10,0){\line(0,1){10}}
\put(20,3){\line(0,1){4}}
\end{picture}
\\
\hline
AII &\begin{picture} (40,-10)
\put(0,0){\line(0,1){10}}
\put(0,0){\line(1,0){20}}
\put(23,0){\line(1,0){4}}
\put(33,0){\line(1,0){4}}
\put(40,0){\line(1,0){10}}
\put(50,0){\line(0,1){10}}
\put(0,10){\line(1,0){20}}
\put(23,10){\line(1,0){4}}
\put(33,10){\line(1,0){4}}
\put(40,10){\line(1,0){10}}
\put(40,3){\line(0,1){4}}
\put(10,0){\line(0,1){10}}
\put(20,3){\line(0,1){4}}
\put(0,-10){\line(0,1){10}}
\put(0,-10){\line(1,0){20}}
\put(23,-10){\line(1,0){4}}
\put(33,-10){\line(1,0){4}}
\put(40,-10){\line(1,0){10}}
\put(50,-10){\line(0,1){10}}
\put(40,-7){\line(0,1){4}}
\put(10,-10){\line(0,1){10}}
\put(20,-7){\line(0,1){4}}
\end{picture}
$\phantom{\begin{array}{c} a\\ b\\ a\end{array}}$
\\
\hline
AIII& $ab....ba, \quad ba....ab, \quad ab....ab, \quad ba....ba$, \\
\hline
BDI& $ab....ba, \qquad ba....ab, \qquad \begin{array}{l}ab....ab\\ ba....ba \end{array}$,\\
\hline
CI& $ab....ab, \qquad ba....ba, \qquad \begin{array}{l}ab....ba\\ ba....ab \end{array}$,\\
\hline
DIII& $\begin{array}{l}ab....ab\\ ab....ab \end{array}, \qquad \begin{array}{l}ba....ba\\ ba....ba \end{array}, \qquad \begin{array}{l}ab....ba\\ ba....ab \end{array},$\\
\hline
CII& $\begin{array}{l}ab....ba\\ ab....ba \end{array}, \qquad \begin{array}{l}ba....ab\\ ba....ab \end{array}, \qquad \begin{array}{l}ab....ab\\ ba....ba \end{array},$\\
\hline
\end{tabular}
\end{center}

\subsection{Cas AI et AII}
La paire $(\g,\kk)$ considérée est de la forme suivante: $(\sln(V),\aut(V,\Phi))$ où $\Phi$ est une forme bilinéaire symétrique (AI) ou anti-symétrique (AII).
On pose $\g'=\gl(V)$ équipée de l'involution définie en \ref{cas} (cas AI et AII). On a alors $\g'=\g+\K \Id_n$ où $\K \Id_n\subset \pp'$.
D'après le lemme \ref{cent}, on a $\g'^\sfr\cong \bigoplus_d \gl_{m_d}$. Et à l'intérieur de ce centralisateur, par la proposition \ref{AI} avec $\wfr=\kk$ et $\eta=-1$, on obtient $\kk^\sfr\cong\bigoplus_d \so_{m_d}$ (AI) ou $\kk^\sfr\cong\bigoplus_d \spn_{m_d}$ (AII).
Pour $d\in \N$ fixé, la paire réductive symétrique associée est donc une paire réductive du même type que celle de l'algèbre de Lie symétrique de départ.
On en déduit que $\delta^e(d)=m_d$ (AI) ou $\delta^e(d)= \frac {m_d} 2$ (AII).
Regardons maintenant dans quels cas $(\g'^\sfr=\kk^\sfr\oplus T_r)$ où $T_r$ est un tore de dimension $r\geqslant 1$, auquel cas $e$ est presque distingué de défaut $r-1$. Il suffit de regarder les sous-algèbres $(\g'_d,\kk_d)$.
\begin{itemize}
\item[AI:] \begin{itemize}
\item Si $m_d=1$, la paire réductive associée est $(T_1,\{0\})$ avec $\delta^e(d)=1$,
\item Si $m_d\geqslant 2$, il existe un élément nilpotent non nul pour la paire $(\gl_{m_d}, \so_{m_d})$. 
Cet élément correspond alors par la théorie des représentations à un élément nilpotent non nul pour la paire $(\gl_n, \so_n)$ 
qui est également nilpotent dans $(\sln_n, \so_n).$
\end{itemize}
\item[AII:] Les éléments $m_d$ sont pairs.
\begin{itemize}
\item Si $m_d=2$, la paire reductive associée est $(\gl_2,\sln_2)$ où $\gl_2\cong\sln_2\oplus T_1$.
\item Si $m_d\geqslant 4$, il existe des éléments nilpotents non nuls pour $(\g'_d,\kk_d)$. Ces éléments correspondent à des éléments nilpotents non nuls de $(\g, \kk)$.
\end{itemize}
\end{itemize}

On en déduit la proposition suivante:
\begin{prop}\label{pdAI}
Les éléments $\pp$-distingués dans le cas AI et AII sont les éléments $\pp$-réguliers. 
Les éléments presque $\pp$-distingués du cas AI (resp. AII) sont ceux qui correspondent à un diagramme de Young ayant des lignes de longueurs distinctes 
(resp. des paires de lignes de longueurs distinctes). Le défaut d'un élément nilpotent est égal respectivement, au nombre de lignes ou de paires de lignes, moins une.
\end{prop}
\begin{proof}
L'assertion sur le défaut provient de la formule $$\delta(e)=\sum_d \delta^e(d)-1$$ de la remarque \ref{delta0}.
Celle sur les éléments $\pp$-distingués provient de la discussion de cas ci-dessus. Enfin, il est connu que les éléments $\pp$-réguliers dans le cas AI (resp. AII) correspondent aux diagrammes de Young constitués d'une seule ligne (resp. paire de lignes).
\end{proof}

\begin{ex}\label{exAI}
Dans le cas AI, les éléments correspondant aux diagrammes suivants sont presque $\pp$-distingués de défaut $1$. En particulier, ils ne sont pas $\pp$-distingués:
$$ \Gamma_1=\begin{picture} (0,0) 
\put(0,-10){\line(0,1){20}}
\put(0,10){\line(1,0){20}}
\put(0,0){\line(1,0){20}}
\put(10,-10){\line(0,1){20}}
\put(20,0){\line(0,1){10}}
\put(0,-10){\line(1,0){10}}
\end{picture} \qquad\quad, \qquad
\Gamma_1'=\begin{picture} (0,0) 
\put(0,-10){\line(0,1){20}}
\put(0,10){\line(1,0){30}}
\put(0,0){\line(1,0){30}}
\put(10,-10){\line(0,1){20}}
\put(20,0){\line(0,1){10}}
\put(30,0){\line(0,1){10}}
\put(0,-10){\line(1,0){10}}
\end{picture}\qquad \qquad, \qquad
\Gamma_2=\begin{picture} (0,0) 
\put(0,-10){\line(0,1){20}}
\put(0,10){\line(1,0){30}}
\put(0,0){\line(1,0){30}}
\put(10,-10){\line(0,1){20}}
\put(20,-10){\line(0,1){20}}
\put(30,0){\line(0,1){10}}
\put(0,-10){\line(1,0){20}}
\end{picture}\qquad \qquad, \qquad
\Gamma_3=\begin{picture} (0,0) 
\put(0,-10){\line(0,1){20}}
\put(0,10){\line(1,0){40}}
\put(0,0){\line(1,0){40}}
\put(10,-10){\line(0,1){20}}
\put(20,-10){\line(0,1){20}}
\put(30,0){\line(0,1){10}}
\put(40,0){\line(0,1){10}}
\put(0,-10){\line(1,0){20}}
\end{picture}\qquad \qquad\quad.
$$
Dans le cas AII, ce sont les mêmes en doublant les lignes, par exemple:
$$ \Gamma_4=\begin{picture} (0,0) 
\put(0,-30){\line(0,1){40}}
\put(0,10){\line(1,0){20}}
\put(0,-10){\line(1,0){20}}
\put(0,0){\line(1,0){20}}
\put(10,-30){\line(0,1){40}}
\put(20,-10){\line(0,1){20}}
\put(0,-20){\line(1,0){10}}
\put(0,-30){\line(1,0){10}}
\end{picture}\qquad \qquad
$$
\end{ex}

\vspace{0.5cm}

\subsection{Cas AIII}
La paire $(\g,\kk)$ considérée est de la forme suivante: $(\sln_n,\sln_p\oplus\sln_q\oplus T_{1})$ où $p+~q~=~n$.
D'après le lemme \ref{cent}, on a $\g^\sfr\cong (\bigoplus_d \gl_{m_d})\cap \sln_n$. Et à l'intérieur de ce centralisateur, par la proposition \ref{AIII}, on obtient $\kk^\sfr\cong\bigoplus_d (\gl_{a_d}\oplus \gl_{b_d})\cap \sln_{m_d}$. 
Pour $d\in \N$ fixé, la paire réductive symétrique associée est donc une paire du même type que celle de l'algèbre de Lie symétrique de départ.
Regardons maintenant dans quels cas $(\g^\sfr=\kk^\sfr\oplus T_r)$ où $T_r$ est un tore de dimension $r\geqslant 0$. 
\begin{itemize}
\item Si $a_db_d=0$, la paire associée est $(\sln_n, \sln_n)$.
\item Si $a_d\neq 0$ et $b_d\neq 0$, il existe un élément nilpotent non nul pour la paire $(\gl_{m_d},(\gl_{a_d}\oplus \gl_{b_d}))$. 
C'est un élément nilpotent pour la paire\\
 $(\sln_n,\sln_p~\oplus~\sln_q~\oplus T_{1})$.
\end{itemize}
On en déduit la proposition suivante:

\begin{prop}\label{pdAIII}
Les éléments $\pp$-distingués dans le cas AIII coïncident avec les éléments presque $\pp$-distingués. Ce sont ceux qui correspondent aux $ab$-diagrammes de Young dont les lignes de longueur fixée commencent par la même lettre.
\end{prop}

\subsection{Cas BDI, CI, DIII et CII}
La paire $(\g,\kk)$ considérée est de la forme suivante: $(\aut(\Phi),\aut(\Phi,J))$ où $\Phi$ est une forme bilinéaire non-dégénérée symétrique 
$(\varepsilon=1)$ (BD) ou antisymétrique $(\varepsilon=-1)$ (C); $J$ preserve $\Phi$ et vérifie $J^2=\xi \Id$ avec $\xi=1$ (BDI, CII) ou $\xi=-1$ (DIII,CI).

D'après la proposition \ref{AI}, on a $\g^\sfr\cong \bigoplus_d \aut(H_d,\Phit_d)$. A l'intérieur de ce centralisateur, par la proposition \ref{BDI}, on obtient $$\kk^s\cong\bigoplus_{d\in \N^*}(\aut(H_d, \Phit_d))^{(\sqrt{-1})^{d-1}J}.$$ 
On pose $\varepsilon_d=1$ si $\Phit_d$ est symétrique et $\varepsilon_{d}=-1$  si $\Phit_d$ est anti-symétrique. 
De même, on définit $\xi_d$ par $((\sqrt{-1})^{d-1}J_{\mid H_d})^2=\xi_d\Id$. 
Par les lemmes \ref{sym} et \ref{J}~(b) où l'on pose $\eta=1$, on a $(\varepsilon_d,\xi_d)=((-1)^{d-1}\varepsilon, (-1)^{d-1}\xi)$. 
Notons que les sous-espaces propres de $(\sqrt{-1})^{d-1}J$ sont de dimensions $a_d$ et $b_d$; si $\xi_d=-1$ alors nécessairement $a_d=b_d=\frac {m_d} 2$. 
Pour $d\in \N$ fixé, la paire réductive symétrique associée est donc donnée par le tableau suivant:
\medskip

$\!\!\!\!\!\!\!\!\!\!\!\!\!\!\!\!\!\!\!\!$\begin{tabular}{|c| c|c|c|c|c|}
\hline
$(\g,\kk)$ &$(\varepsilon, \xi)$ &  \multicolumn{2}{|c|}{$d$ pair}&  \multicolumn{2}{|c|}{$d$ impair}\\
\hline
& & $(\varepsilon_d, \xi_d)$ & $(\g_d,\kk_d)$  & $(\varepsilon_d, \xi_d)$ & $(\g_d,\kk_d)$\\
\hline
$(\so_n,\so_p\times\so_q)$ (BDI)& (1,1) & (-1,-1) &$(\spn_{m_d}, \gl_{\frac {m_d} 2})$ & (1,1) & $(\so_{m_d},\so_{a_d}\times\so_{b_d})$ \\
\hline
$(\spn_n,\gl_{\frac n 2})$ (CI) & (-1,-1)& (1,1) &$(\so_{m_d},\so_{a_d}\times\so_{b_d})$ &  (-1,-1) &$(\spn_{m_d}, \gl_{\frac {m_d} 2})$ \\
\hline
$(\so_n,\gl_{\frac n 2})$ (DIII) & (1,-1) & (-1,1) &$(\spn_{m_d}, \spn_{a_d}\times \spn_{b_d})$ & (1,-1) &$(\so_{m_d},\gl_{\frac {m_d} 2})$\\
\hline
$(\spn_n,\spn_p\times \spn_q)$ (CII) & (-1,1) & (1,-1) &$(\so_{m_d},\gl_{\frac {m_d} 2})$ & (-1,1) & $(\spn_{m_d}, \spn_{a_d}\times \spn_{b_d})$\\
\hline
\end{tabular}
\medskip

Regardons maintenant dans quels cas $(\g^\sfr=\kk^\sfr\oplus T_r)$ où $T_r$ est un tore de dimension $r\geqslant 0$.
\begin{itemize}
\item[(a)] (BDI) Si $a_db_d=0$, alors $(\g_d, \kk_d)=(\so_{m_d}, \so_{m_d})$.
\item[(b)] (BDI) Si $a_db_d=1$, on a $(\g_d, \kk_d)=(T_1, \{0\})$.
\item[(c)] (BDI) Si $a_db_d\geqslant 2$, $\pp_d$ contient des éléments nilpotents non nuls.
\item[(d)] (CI) Si $m_d\geqslant2$ (i.e. $\neq 0$), $\pp_d$ contient des éléments nilpotents non nuls.
\item[(e)] (CII) Si $a_db_d=0$ alors $(\g_d, \kk_d)=(\spn_{m_d}, \spn_{m_d})$.
\item[(f)] (CII) Si $a_d,b_d\geqslant1$, $\pp_d$ contient des éléments nilpotents non nuls.
\item[(g)] (DIII) Si $m_d=2$, on a $(\g_d, \kk_d)=(T_1, T_1)$.
\item[(h)] (DIII) Si $m_d\geqslant 4$, $\pp_d$ contient des éléments nilpotents non nuls.
\end{itemize}

On en déduit les deux propositions suivantes:
\begin{prop}\label{pdBDI}
Les éléments $\pp$-distingués du cas BDI (resp. CI) sont ceux qui ont un $ab$-diagramme n'ayant pas de lignes de longueur paire (resp. impaire) et tels que les lignes de longueur impaire (resp.paire) débutent toutes par la même lettre. Pour obtenir les éléments presque $\pp$-distingués, on autorise de plus, par longueur de ligne impaire (resp. paire) une paire de ligne, débutant l'une par $a$ et l'autre par $b$. Le nombre de tels couples est égal au défaut de notre élément.
\end{prop}

\begin{proof}
Dans le cas (a), on a $\delta^e(d)=0$. Dans le cas (b), on a $\delta^e(d)=\dim \pp^{\sfr}_d=1$. Dans les cas (c) et (d), on a des éléments nilpotents non nuls 
pour la paire $(\g^{\sfr},\kk^{\sfr})$, ce qui empêche l'élément d'être presque $\pp$-distingué. On en déduit la proposition grâce au tableau précédent.
\end{proof}

\begin{prop}\label{pdDIII}
Les éléments $\pp$-distingués du cas CII (resp. DIII) sont ceux qui ont un $ab$-diagramme tels qu'il existe, par longueur de ligne paire (resp. impaire), 
une ou zéro paire de ligne (débutant nécessairement par $a$ et par $b$) et tels que les lignes de longueur impaire (resp. paire) débutent toutes par la même lettre. 
Par ailleurs, les éléments presque $\pp$-distingués sont $\pp$-distingués.
\end{prop}

\begin{proof}
Dans les cas (e) et (g), on a $\delta^e(d)=0$. Dans les cas (f) et (h), on a des éléments nilpotents non nuls pour la paire $(\g^{\sfr},\kk^{\sfr})$, 
ce qui empêche l'élément d'être presque $\pp$-distingué. On en déduit la proposition grâce au tableau précédent.
\end{proof}

\begin{ex}\label{exBDI}
Dans le cas BDI, les éléments correspondant aux diagrammes suivants sont presque $\pp$-distingués de défaut $1$. En particulier, ils ne sont pas distingués:
$$ \Gamma_5=\begin{array}{l} aba \\a\\ b\end{array},\qquad \Gamma_6=\begin{array}{l} ababa \\aba\\ bab\\ b \end{array},\qquad \Gamma_7=\begin{array}{l} ababa \\aba\\ bab\\ a\end{array}.$$
Le cas CI est très similaire, des exemples de tels éléments sont les suivants:
$$\Gamma_8=\begin{array}{l} baba \\ba\\ ab\end{array},\qquad \Gamma_9=\begin{array}{l} bababa \\baba\\ abab\\ ab \end{array},\qquad \Gamma_{10}=\begin{array}{l} bababa \\baba\\ abab\\ ba\end{array}.$$ 
\end{ex}

\subsection{Bilan Provisoire}\label{bilan1}
A partir des examens de $(\g'_d,\kk_d)$ dans certains cas classiques simples effectués lors de cette section, on peut résumer le défaut de longueur $d$ associé à 
$e$ dans le tableau ci-dessous.
\begin{center}
\begin{tabular}{|c|c|c|}
\hline
Type & \multicolumn{2}{|c|}{$\delta^e(d)$ pour $d\in \N^*$}\\
\hline
AI & \multicolumn{2}{|c|}{$m_d$}\\
\hline
AII & \multicolumn{2}{|c|}{$\frac {m_d} 2$}\\
\hline
& $d$ pair & $d$ impair \\
\hline
(BDI)& $\frac{m_d} 2$  & $\min\{a_d,b_d\}$ \\
\hline
(CI) & $\min\{a_d,b_d\}$ & $\frac{m_d} 2$\\
\hline
\end{tabular}
\end{center}

Par ailleurs, en combinant les résultats \ref{pdAIII} et \ref{pdDIII} avec la proposition \ref{pd}, on obtient:
\begin{cor} La conjecture \ref{conj} est vraie dans les cas AIII, DIII et CII. \end{cor}

\section{Réduction d'orbites presque $\pp$-distinguées}
\label{reduction}

Nous disposons d'un autre résultat pour écarter le fait que certains éléments presque $\pp$-distingués puissent engendrer une composante irréductible étrange.

\begin{prop}\label{xx}
\begin{description}
\item (a) Soit $e_1, e_2\in \pp$ deux éléments nilpotents avec $\Od(e_1)\subset\overline{\Od(e_2)}$ alors $\delta(e_1)-\delta(e_2)\leqslant\dim \Od(e_2)-\dim \Od(e_1)=\dim \pp^{e_{1}}-\dim \pp^{e_{2}}$.
\item (b) Si il y a égalité et $e_1$ est presque $\pp$-distingué, alors $\CC(e_1)\subset\CC(e_2)$.
\end{description} 
\end{prop}
\begin{proof}
 (a) Rappelons que l'on a  une application $\pr_1:\CC(e_2)\rightarrow\overline{\Od(e_2)}$.
 Comme $(g.e_2,g.e_2)\in \CC(e_2)$ pour tout $g\in K$ et comme $e_1\in \overline{\Od(e_2)}$, on a $(e_1,e_1)\in \CC(e_2)$ et $e_1$ appartient à l'image de $\pr_1$. Maintenant, pour tout $x\in \Od(e_2)$, on a $\pr_1^{-1}(x)=(x,\pp^x\cap\NN)$. On en déduit donc que $\pr_1^{-1}(e_1)\subseteq (e_1,\pp^{e_1}\cap~\NN)$ est un sous ensemble de dimension au moins égale à $\dim \pp^{e_2}\cap\NN$.
On en déduit que 
$\dim\pp^{e_2}-\delta(e_2)=\dim\pp^{e_2}\cap\NN\leqslant\dim\pp^{e_1}\cap\NN=\dim\pp^{e_1}-\delta(e_1).$
L'inégalité de (a) s'en suit facilement.

 (b) Si les hypothèses de (b) sont vérifiées, alors $\pp^{e_1}\cap\NN=\bigoplus_{i\geqslant1}\pp(e_{1},i)$ est un ensemble irréductible contenant une fibre de $\pr_1$
 de même dimension, d'où $(e_1,\pp^{e_1}\cap\NN)\subset\CC(e_2)$ et $\CC(e_1)\subset\CC(e_2)$.
\end{proof}

\begin{defi}
Si $e_1$ est presque $\pp$-distingué et $e_2$ nilpotent quelconque vérifiant l'égalité $\delta(e_1)-\delta(e_2)=\dim \pp^{e_{1}}-\dim \pp^{e_{2}}$, 
on dit que $\Od(e_1)$ se réduit en $\Od(e_2)$ ou encore que $\Od(e_2)$ est une réduction d'ordre $\delta(e_1)-\delta(e_2)$ de $\Od(e_1)$.
\item Une telle réduction est dite minimale si il n'existe pas d'orbite $\Od$ telle que $\overline{\Od(e_1)}\subsetneq\overline{\Od}\subsetneq\overline{\Od(e_2)}$.
\end{defi}

\begin{lm}
Si une orbite $\Od(e_1)$ se réduit en une orbite $\Od(e_2)$ alors toute orbite $\Od(e_3)$ telle que $\overline{\Od(e_1)}\subsetneq\overline{\Od(e_3)}\subsetneq\overline{\Od(e_2)}$, est une réduction de $\Od(e_1)$.
\end{lm}
\begin{proof}
Par la proposition \ref{xx} (a), on a 
\begin{eqnarray*}\delta(e_1)-\delta(e_3) &=&  \big(\delta(e_1)-\delta(e_2)\big)-\big(\delta(e_3)-\delta(e_2)\big)\\
&\geqslant& \big(\dim\pp^{e_1}-\dim\pp^{e_2}\big)-\big(\dim\pp^{e_3}-\dim\pp^{e_2}\big)\\
&\geqslant& \dim\pp^{e_1}-\dim\pp^{e_3}.
\end{eqnarray*}
L'inégalité inverse découle de \ref{xx} (a).
\end{proof}

En vertu de la proposition et du lemme précédents, pour montrer à l'aide d'une réduction qu'un élément $e_1$ presque $\pp$-distingué 
n'engendre pas de composante irréductible étrange, il suffit de savoir qu'il existe une réduction minimale de $e_1$.

Dans chaque cas classique simple non encore résolu, on introduit un ordre sur les ($ab$-)diagrammes pour transférer la relation d'inclusion des orbites sur les diagrammes 
(cf. \cite[1.9]{Ot1}, \cite[1.4]{Ot2}). C'est notamment l'objet des définitions suivantes.
\begin{defi} Soit $\Gamma$ un ($ab$-)diagramme de Young. 
On appelle motif de $\Gamma$ un sous-diagramme obtenu à partir de $\Gamma$ en conservant uniquement les lignes de longueur $d\in I$ où $I$ est un intervalle de $\N$.
\item On note $\Gamma'$ le diagramme obtenu en enlevant la première colonne de $\Gamma$. 
On définit par récurrence $\Gamma^{(k)}=(\Gamma^{(k-1)})'$ pour $k\in \N^*$ où $\Gamma^{(0)}=\Gamma$. 
\item Soit $(\Gamma_{1},\Gamma_{2})$ un couple de ($ab$-)diagrammes de Young. On définit une relation d'ordre partiel comme suit. 
On dit que $\Gamma_{1}\leqslant\Gamma_{2}$ ou que $\Gamma_{2}$ est une dégénérescence de $\Gamma_{1}$ si:\begin{itemize}
   \item pour tout $k \in \N^*$, le nombre de cases de $\Gamma_1^{(k)}$ est inférieur à celui de $\Gamma_2^{(k)}$, 
   lorsque $\Gamma_{1},\Gamma_{2}$ sont des diagrammes de Young;
   \item pour tout $k \in \N^*$, le nombre de $a$ et le nombre de $b$ de $\Gamma_1^{(k)}$ sont inférieurs respectivements à ceux de $\Gamma_2^{(k)}$, 
   lorsque $\Gamma_{1},\Gamma_{2}$ sont des $ab$-diagrammes de Young.\end{itemize}
\item La dégénérescence $\Gamma_{1}\lnq\Gamma_{2}$ est dite minimale s'il n'existe pas d'($ab$-)diagramme $\Gamma$ tel que $\Gamma_{1}\lnq\Gamma\lnq\Gamma_{2}$.
\item On appelle ligne commune à $\Gamma_{1}$ et $\Gamma_{2}$\begin{itemize}
          \item une ligne de même longueur si $\Gamma_{1}$ et $\Gamma_{2}$ sont des diagrammes de Young;
          \item une ligne de même longueur débutant par la même lettre si $\Gamma_{1},\Gamma_{2}$ sont des $ab$-diagrammes de Young. \end{itemize}
\item 
On définit $\overline{\Gamma_1}$ et $\overline{\Gamma_2}$, en enlevant toutes les lignes communes à $\Gamma_1$ et $\Gamma_2$.
\end{defi}

Rappelons que si $\Od_1,\Od_2$ sont deux orbites ayant pour ($ab$-)diagrammes respectifs $\Gamma_{1},\Gamma_{2}$, alors $\Od_{1}\subseteq\overline{\Od_{2}}$ implique $\Gamma_{1}\leqslant\Gamma_{2}$. 
Réciproquement, si deux
\mbox{($ab$-)}diagrammes vérifient $\Gamma_{1}\leqslant\Gamma_{2}$, et si $\Od_1$ a pour diagramme $\Gamma_{1}$, alors il existe $\Od_2$ correspondant à $\Gamma_{2}$ 
telle que $\Od_1\subseteq\overline{\Od_2}$ (cf. \cite[Théorème 3]{Ot1}, \cite[Théorème 1]{Ot2}).

Notons que le fait qu'une orbite se réduise en une autre ne fait intervenir que des invariants propres aux ($ab$-)diagrammes de Young: 
défaut, dimension de centralisateur, inclusion. On transfert alors également les notions de défaut et de dimension de centralisateur aux ($ab$-)diagrammes 
que l'on note repectivement $\delta(\Gamma)$ et $\dim \pp^{\Gamma}$.
\begin{defi} Soient $\Gamma_1,\Gamma_2$ deux ($ab$-)diagrammes de Young, on pose 
 $\Delta=\Delta(\Gamma_{1},\Gamma_{2})=\delta(\Gamma_1)-\delta(\Gamma_2)$ et $s=s(\Gamma_{1},\Gamma_{2})=\dim\pp^{\Gamma_1}-\dim \pp^{\Gamma_2}$.
\item On dit que $\Gamma_1$ se réduit en $\Gamma_2$, ou que $\Gamma_{2}$ est une réduction de $\Gamma_{1}$, si $\Gamma_1\lnq\Gamma_2$ et $s=\Delta$.
\end{defi}

\begin{Rq} Les notions de défaut, de dimension de centralisateur et de réduction pour les ($ab$-)diagrammes dépendent du cas simple considéré (AI, AII, BDI ou CI).
\item Soit $\Od_{1}$ une orbite et $\Gamma_{1}$ son ($ab$-)diagramme. L'existence d'une réduction pour $\Od_{1}$ équivaut à celle d'une réduction pour $\Gamma_{1}$.
\item La notion de dégénérescence minimale équivaut à celle d'\og{adjacent degeneration}\fg introduite dans \cite[(1.4)]{Ot2} et \cite[(2.4)]{Ot1}.
\end{Rq}

\subsection{Les cas AI et AII}
Regardons le cas le plus simple: le cas AI.
On fixe un élément $e_1$ presque $\pp$-distingué, et on note son diagramme de Young $\Gamma_1$. 
Par la proposition \ref{pdAI}, les lignes de $\Gamma_1$ sont de longueurs distinctes. 
D'après \cite[Lemme 5]{Ot2}, la seule dégénéréscence minimale $\Gamma_{2}$ de $\Gamma_{1}$ susceptible d'être une réduction doit fournir 
\medskip
$$
\overline{\Gamma_1}:
\begin{picture} (40,-10)
\put(0,0){\line(0,1){10}}
\put(0,0){\line(1,0){10}}
\put(13,0){\line(1,0){4}}
\put(23,0){\line(1,0){4}}
\put(30,0){\line(1,0){10}}
\put(40,0){\line(0,1){10}}
\put(0,10){\line(1,0){10}}
\put(13,10){\line(1,0){4}}
\put(23,10){\line(1,0){4}}
\put(30,10){\line(1,0){10}}
\put(30,3){\line(0,1){4}}
\put(10,3){\line(0,1){4}}
\put(0,-10){\line(0,1){10}}
\put(0,-10){\line(1,0){10}}
\put(13,-10){\line(1,0){4}}
\put(23,-10){\line(1,0){4}}
\put(30,-10){\line(1,0){10}}
\put(40,-10){\line(0,1){10}}
\put(30,-7){\line(0,1){4}}
\put(10,-7){\line(0,1){4}}
\put(40,10){\line(0,-1){20}}
\put(40,10){\line(1,0){20}}
\put(40,0){\line(1,0){20}}
\put(70,0){\line(1,0){10}}
\put(63,0){\line(1,0){4}}
\put(70,10){\line(1,0){10}}
\put(63,10){\line(1,0){4}}
\put(40,-10){\line(1,0){10}}
\put(50,10){\line(0,-1){20}}
\put(80,10){\line(0,-1){10}}
\put(70,7){\line(0,-1){4}}
\put(42,-7){x}
\end{picture}
\hspace{3cm}
\overline{\Gamma_2}:
\begin{picture} (40,-10)
\put(0,0){\line(0,1){10}}
\put(0,0){\line(1,0){10}}
\put(13,0){\line(1,0){4}}
\put(23,0){\line(1,0){4}}
\put(30,0){\line(1,0){10}}
\put(40,-10){\line(0,1){20}}
\put(0,10){\line(1,0){10}}
\put(13,10){\line(1,0){4}}
\put(23,10){\line(1,0){4}}
\put(30,10){\line(1,0){10}}
\put(30,3){\line(0,1){4}}
\put(10,3){\line(0,1){4}}
\put(0,-10){\line(0,1){10}}
\put(0,-10){\line(1,0){10}}
\put(13,-10){\line(1,0){4}}
\put(23,-10){\line(1,0){4}}
\put(30,-10){\line(1,0){10}}
\put(30,-7){\line(0,1){4}}
\put(10,-7){\line(0,1){4}}
\put(40,10){\line(1,0){20}}
\put(40,0){\line(1,0){20}}
\put(70,0){\line(1,0){10}}
\put(63,0){\line(1,0){4}}
\put(70,10){\line(1,0){10}}
\put(63,10){\line(1,0){4}}
\put(50,10){\line(0,-1){10}}
\put(80,10){\line(0,-1){10}}
\put(70,7){\line(0,-1){4}}
\put(80,0){\line(1,0){10}}
\put(80,10){\line(1,0){10}}
\put(90,0){\line(0,1){10}}
\put(82,3){x}
\end{picture}
$$
\medskip


C'est à dire que si l'on note $\Gamma_i=(p^i_1,p^i_2,\dots, p^i_n)$ où $(p^i_j)_j$ est une suite décroissante d'éléments de somme $n$,
alors il existe $j_0\in[\![1,n]\!]$ tel que $$p_{j_0}^2=p_{j_0}^1+1,\qquad p_{j_0+1}^2=p_{j_0+1}^1-1,\qquad p^2_j=p^1_j \mbox{ sinon. }$$
La dimension du centralisateur dans $\pp$ d'un élément correspondant à un tel diagramme de Young est donnée par (cf. \cite[3.1]{Se}):
\begin{equation} \label{dimAI} \dim \pp^{\Gamma_i}=(\sum_{j\in \N^*} jp_j^i)-1. \end{equation}

On peut donc calculer la différence $s$:
$$s=\dim \pp^{\Gamma_1}-\dim \pp^{\Gamma_2}=j_0p^1_{j_0}+(j_0+1)p^1_{j_0+1}-j_0(p^1_{j_0}+1)-(j_0+1)(p^1_{j_0+1}-1)=1.$$
Pour obtenir une réduction, on doit donc avoir $\Delta=1$, c'est à dire $\delta(\Gamma_1)=\delta(\Gamma_2)+1$ ce qui se traduit par le fait que $\Gamma_2$
possède une ligne de moins que $\Gamma_1$ (cf. proposition \ref{pdAI}). Ceci n'est possible que si $\Gamma_1$ possède une ligne (nécessairement unique) de longueur 1.
Ainsi, si l'on reprend les diagrammes de l'exemple \ref{exAI}, on voit que $\Gamma_1$ et $\Gamma_1'$ possèdent chacun une réduction, 
mais que ce n'est le cas ni pour $\Gamma_2$, ni pour $\Gamma_3$. Au passage, ceci permet de prouver la conjecture \ref{conj} dans le cas AI de rang 2 et 3, 
car dans ces cas $\Gamma_1$ et $\Gamma_1'$ représentent les seuls éléments presque $\pp$-distingués non $\pp$-distingués. 
\newline

Dans le cas AII, par la même méthode, nous pouvons montrer que la seule dégénérescence minimale possible pour un élément presque $\pp$-distingué est obtenue de celle du cas AI en doublant les lignes. On a alors $s=4$ tandis que $\Delta=1$, on n'a donc pas de réduction d'élément presque $\pp$-distingué.

\begin{cons} La conjecture \ref{conj} est démontrée dans le cas AI jusqu'au rang 3 et dans le cas AII de rang 1.\end{cons}

\subsection{Les cas BDI et CI}

Heureusement, dans les cas BDI et CI, la méthode de réduction permet d'obtenir plus de résultats.
On traite tout d'abord le cas BDI, le cas CI lui étant très similaire. 

Introduisons quelques notions pour calculer les entiers $s(\Gamma_{1},\Gamma_{2})$. On fixe un élément $e$ nilpotent quelconque (non nul) 
dans une algèbre symétrique de type BDI, et $\Gamma$ désigne son $ab$-diagramme de Young.

\begin{defi}
Si $\Gamma$ ne comporte que des lignes de longueur impaire on définit $k^{\Gamma}_j$ comme étant la longueur de la $(2j+1)$-ème colonne de $\Gamma$.
\end{defi}

On rappelle que le défaut de la longueur $d$ est donné par $\delta^{\Gamma}(d)=\min(a_d, b_d)$ et que $\delta(\Gamma)=\sum_d \delta^{\Gamma}(d)$ 
(cf. section \ref{bilan1} et remarque \ref{delta0}).

\begin{lm}\label{mmm}
Si $h$ est un élément à valeurs propres entières de $\so(V)$ et $V(i)$ le $h$-module de poids $i\in \Z$, alors 
$\g^h\cong\so(V(0))\oplus\bigoplus_{i\gnq0} \gl(V(i))$.
\end{lm}
\begin{proof}
On note $\Phi$ la forme bilinéaire symétrique canonique définissant $\so(V)$. Soit $v\in V(i)$ et $w\in V(j)$. On a $\Phi(h.v,w)+\Phi(v,h.w)=0$, d'où $\Phi(v,w)=0$ ou $i+j=0$. On a donc une décomposition orthogonale $h$-invariante $V=V(0)\oplus\bigoplus_{i\in\N^*} (V(i)\oplus V(-i))$. On voit que $(\g^h)_{\mid V(i)\oplus V(-i)}\cong \gl(V(i))$. D'où le résultat.
\end{proof}

\begin{lm}\label{dimBDI}
Si $\Gamma$ ne comporte que des lignes de longueur impaire, on a $$\dim \pp^{\Gamma}=C+\frac{k_0^{\Gamma}(k_0^{\Gamma}-1)} {4} +\frac12\sum_{j\gnq0} (k_j^{\Gamma})^2$$ où $C$ est une constante qui ne dépend que de l'algèbre symétrique considérée.
\end{lm}
\begin{proof}
On inclut $e$ dans un $\Striplet$ normal $(e,h,f)$. L'élément $e$ est pair donc $\dim \pp^e=(\dim \pp-\frac{1}{2}\dim\g)+ \frac{1}{2}\dim\g^h.$ 
Comme $\Gamma$ ne contient que des lignes de longueur impaire, on peut décomposer $V$ en $\bigoplus_{j\in \Z} V(2j)$. 
On voit alors par la théorie des représentations de $\sld$ que $k^{\Gamma}_j=\dim V(2j)$ pour $j\in \N$.
On obtient ensuite par le lemme \ref{mmm} que $\frac12\dim \g^h=\frac{k_0^{\Gamma}(k_0^{\Gamma}-1)} {4} +\frac12\sum_{j\gnq0} (k_j^{\Gamma})^2.$
\end{proof}

On fixe maintenant un élément presque $\pp$-distingué $e_1$ et on note $\Gamma_1$ son $ab$-diagramme de Young. On énumère les dégénérescences minimales possibles 
(cf. \cite[Table~V cas BDI]{Ot1}). La première colonne désigne le numéro de la dégénérescence minimale $\Gamma_1\lnq \Gamma_2$ 
telle qu'elle est indiquée dans \cite{Ot1}. Les deux colonnes suivantes donnent $\overline{\Gamma_1}$ et $\overline{\Gamma_2}$. 
La quatrième donne l'entier $s=\dim\pp^{\Gamma_1}-\dim \pp^{\Gamma_2}$ et la dernière indique les restrictions sur les entiers $p$ et $q$.

Notons que le cas (1) de \cite{Ot1} n'apparaît pas car $\Gamma_1$ ne peut alors pas correspondre à une orbite presque $\pp$-distingué. 
Les cas (6) et (7) n'apparaissent pas car, pour ces dégénérescences, on a $\Delta\leqslant 0$. 
L'entier $s$ a été calculé à partir du lemme \ref{dimBDI}, sauf dans certains cas dont on verra plus tard qu'ils ne peuvent donner lieu à de nouvelles réductions.
\newline

\begin{tabular}{|p{1.5cm}|c|c|c|c|}
\hline
numéro de \cite{Ot1} & \rule{0pt}{10pt} $\overline{\Gamma_1}$ & $\overline{\Gamma_2}$ &$s$& précisions\\ 
\hline
(2) & \begin{tabular}{l}$\overbrace{ab.....ba}^{2p+1}$\\ $\underbrace{ab...ba}_{2q+1}$\end{tabular} & \begin{tabular}{l}$\overbrace{ab.....ba}^{2p+3}$\\ $\underbrace{ab...ba}_{2q-1}$\end{tabular} & $k^{\Gamma_{1}}_q-k^{\Gamma_{1}}_{p+1}-1$ & $p\geqslant q\geqslant 1$\\
\hline
(3) & \begin{tabular}{l}$\overbrace{ab.....ba}^{2p+1}$\\ $\underbrace{ba...ab}_{2q+1}$\end{tabular} & \begin{tabular}{l}$\overbrace{ab.....ba}^{2p+3}$\\ $\underbrace{ba...ab}_{2q-1}$\end{tabular} &$k^{\Gamma_{1}}_q-k^{\Gamma_{1}}_{p+1}-1$& $p\geqslant q\geqslant 1$\\ 
\hline
(4) & \begin{tabular}{l}$\overbrace{ab.....ba}^{2p+1}$\\ $ba.....ab$\\ $\underbrace{ab...ba}_{2q+1}$\end{tabular} & \begin{tabular}{l}$\overbrace{ab.....ab}^{2p+2}$\\$ba.....ba$\\ $\underbrace{ab...ba}_{2q-1}$\end{tabular} & - - - - & $p\geqslant q\geqslant 1$\\
\hline
(5) & \begin{tabular}{l}$\overbrace{ab.....ba}^{2p+1}$\\ $ab...ba$\\ $\underbrace{ba...ab}_{2q+1}$\end{tabular} & \begin{tabular}{l}$\overbrace{ab.....ba}^{2p+3}$\\$ab...ab$\\ $\underbrace{ba...ba}_{2q}$\end{tabular} & \begin{tabular}{c}$k^{\Gamma_{1}}_0-k^{\Gamma_{1}}_{p+1}-2$\\- - - - \end{tabular} & \begin{tabular}{c}$p\geqslant q=0$\\ $p\geqslant q\geqslant 1$ \end{tabular}\\
\hline
(8) & \begin{tabular}{l}$\overbrace{ba.....ab}^{2p+1}$\\ $\underbrace{ab...ba}_{2q+1}$\end{tabular} & \begin{tabular}{l}$\overbrace{ab.....ba}^{2p+3}$\\ $\underbrace{ba...ab}_{2q-1}$\end{tabular} &$k^{\Gamma_{1}}_q-k^{\Gamma_{1}}_{p+1}-1$& $p\geqslant q\geqslant1$\\
\hline
(9) & \begin{tabular}{l}$\overbrace{ba.....ab}^{2p+1}$\\ $ab...ba$\\ $\underbrace{ab...ba}_{2q+1}$\end{tabular} & \begin{tabular}{l}$\overbrace{ba.....ab}^{2p+3}$\\$ba...ba$\\ $\underbrace{ab...ab}_{2q}$\end{tabular} & \begin{tabular}{c}$k^{\Gamma_{1}}_0-k^{\Gamma_{1}}_{p+1}-2$\\- - - - \end{tabular} & \begin{tabular}{c}$p\geqslant q=0$\\ $p\geqslant q\geqslant 1$ \end{tabular}\\
\hline
(10) & \begin{tabular}{l}$\overbrace{ba.....ab}^{2p+1}$\\ $ba.....ab$\\ $\underbrace{ab...ba}_{2q+1}$\end{tabular} & \begin{tabular}{l}$\overbrace{ba.....ba}^{2p+2}$\\$ab.....ab$\\ $\underbrace{ab...ba}_{2q-1}$\end{tabular} & - - - - & $p\geqslant q\geqslant 1$\\
\hline
\end{tabular}
\newline

Pour une dégénérescence minimale donnée, on note $\Delta(d)=\delta^{\Gamma_1}(d)-\delta^{\Gamma_2}(d)$.
Notons que $\Delta(d)\in\{0,1,-1\}$ par le tableau de la section \ref{bilan1}.
\newline

Considérons la dégénérescence minimale (3). D'après le lemme \ref{dimBDI}, on a \begin{equation} \label{s}s=k_q-k_{p+1}-1.\end{equation}
On voit que les réductions éventuelles de défaut ne peuvent se produire que pour les longueurs $2q+1$ et $2p+1$. Les augmentations éventuelles de défaut, pour les longueurs $2p+3$ et $2q-1$. 
\newline

On va tout d'abord se placer dans le cas $p=q$. On a alors $a_{2p+1}b_{2p+1}\neq0$, or $e_1$ est presque distingué donc $a_{2p+1}b_{2p+1}=1$. On en déduit $k_q=k_p=k_{p+1}+2$ et $s=1$ d'après (\ref{s}). Or $\Delta(2p+1)=1$. Pour obtenir une réduction, on doit nécéssairement avoir $\Delta(2p+3)=\Delta(2p-1)=0$. Voyons ce que cela donne dans les différents cas: 
\begin{itemize} 
\item Si $a_{2p+3}b_{2p+3}=1$ ou $b_{2p+3}=0$ (resp. $a_{2p-1}b_{2p-1}=1$ ou $a_{2p-1}=0$), alors $\Delta(2p+3)=0$ (resp. $\Delta(2p-1)=0$). 
\item Le seul cas pour lequel (3) ne donne pas de réduction est celui où $b_{2p+3}=m_{2p+3}\neq0$ ou $a_{2p-1}=m_{2p-1}\neq0$.
\item On obtient un résultat analogue si l'on considère la réduction (3') où les rôles de $a$ et $b$ sont inversés. On voit alors que le seul cas qui ne peut pas être réduit 
par (3) ou (3') est (à permutation de $a$ et $b$ près) le cas où $a_{2p+3}=m_{2p+3}\neq0$ et $a_{2p-1}=m_{2p-1}\neq0$. Le plus petit exemple est le suivant:
$$\begin{array}{l}
ababa\\aba\\bab\\a
\end{array}.$$
\end{itemize}

Notons que l'on a ainsi montré que si des longueurs de lignes consécutives impaires avaient un défaut, on pouvait réduire l'orbite $\Od(e_1)$. 
Nous supposerons maintenant dans la suite que $\Gamma_1$ ne possède pas de motif
$$\begin{array}{l}\overbrace{ab.....ba}^{2p+3}\\ ba.....ab\\ ab...ba \\ \underbrace{ba...ab}_{2p+1}\end{array}.$$

Nous allons aussi montrer que (3) pour $p\neq q$ ne peut pas apporter de réduction. 
En effet, soit $\Gamma_2$ le diagramme obtenu grâce à une dégénéréscence (3) où l'on suppose $p\neq q$. 
Avec les notations précédentes, il est facile de voir que si $\Delta(2p+1)=1$ (resp. $\Delta(2q+1)=1$) alors $k_{p+1}=k_p-2$ (resp. $k_q=k_p+2$). 
On en déduit alors que 
$$k_q\geqslant k_p+1+\Delta(2q+1), \qquad k_{p+1}\leqslant k_p-1-\Delta(2p+1),$$
$$s= k_q-k_{p+1}-1\geqslant (1+\Delta(2q+1))+(1+\Delta(2p+1))-1 \geqslant \Delta+1.$$
On ne peut donc pas espérer de réduction. Le cas de l'opération (2) se traite de la même manière.
\newline

Regardons maintenant ce qui se passe dans le cas de la dégénérescence minimale (5) pour $q=0$ et $p$ quelconque. 
On remarque que $a_1b_1=1$ car $e_1$ est presque $\pp$-distingué. D'après le lemme \ref{dimBDI}, on a $s=k_0-k_{p+1}-2$. 
Observons que $k_0=k_p+2$ et $\Delta(1)=1$. Notons que le cas où $\delta^{\Gamma_1}(2p+1)=1$ a déjà été montré comme réductible 
(suivant que $p=1$ auquel cas on a deux défauts consécutifs, ou $p\gnq 1$ auquel cas $m_{2p-1}=0$). 
On supposera donc que $\delta^{\Gamma_1}(2p+1)=0$ et on obtient une réduction par (5) si et seulement si $k_{p+1}=k_p-1$ c'est à dire 
$m_{2p+1}=1$. Le seul cas que l'on arrive pas à traiter est le cas où $m_{2p+1}\gnq1$ dont le plus petit exemple est:  
$$\begin{array}{l} aba\\aba\\a\\b
\end{array}.$$

On peut montrer que les autres dégénérescences minimales ne mènent pas à des réductions des exemples que l'on a cité.
Donnons ici des exemples de réduction. On considère les diagrammes $\Gamma_5$ et $\Gamma_6$ de l'exemple \ref{exBDI} et on indique de l'autre coté de la flèche un diagramme qui les réduit.
$$\begin{array}{l} aba \\a\\ b\end{array}\rightarrow ababa; \qquad \begin{array}{l} ababa \\aba\\ bab\\ b \end{array}\rightarrow\begin{array}{l} ababa \\ababa\\ b\\ b \end{array}.$$
En conclusion pour le cas BDI: on peut réduire $e_1$ si et seulement si les lignes de son $ab$-diagramme comportant un défaut ne sont pas toutes dans des motifs du type (à permutation de $a$ et $b$ près):

\begin{center}
\begin{tabular}{l}$\overbrace{\vdots{\phantom{b.........ba}}}^{2p+5}$  \\ $ab.........ba$\\ $ab......ba$\\ $ba......ab$ \\ $ab...ba$\\ $\underbrace{\underbrace{\vdots\phantom{b...ba}}_{2p+1}{\phantom{...}}}_{2p+3}$ \end{tabular} $\qquad$ ou $\qquad$
\begin{tabular}{l}$\overbrace{\vdots{\phantom{b...ba}}}^{2p+3}$ \\ $ab...ba$\\ $ab...ba$\\ $a$ \\ $b$\end{tabular} $\qquad$
pour $p\geqslant 0$.  
\end{center}

Dans le cas CI, les dégénérescences minimales possibles sont indexées dans \cite{Ot1} de (1) à (10) de façon similaire au cas BDI. 
La dégénérescence (5) ne donne plus de réduction, seule, la (3) en donne une. Il existe un unique motif non réductible (à permutation de $a$ et $b$ près):

\begin{center}
\begin{tabular}{l}$\overbrace{\vdots{\phantom{b.........ab}}}^{2p+6}$  \\ $ab.........ab$\\ $ab......ab$\\ $ba......ba$ \\ $ab...ab$\\ $\underbrace{\underbrace{\vdots\phantom{b...ba}}_{2p+2}{\phantom{...}}}_{2p+4}$ \end{tabular} $\qquad$ 
pour $p\geqslant 0$.  
\end{center}

Indiquons, comme ci-dessus, comment les diagrammes $\Gamma_{8}$ et $\Gamma_{9}$ de l'exemple \ref{exBDI} se réduisent:
$$\begin{array}{l} baba \\ba\\ ab\end{array}\rightarrow \begin{array}{l}bababa\\ ab\end{array}; \qquad \begin{array}{l} bababa \\baba\\ abab\\ ab \end{array}\rightarrow\begin{array}{l} bababa \\bababa\\ ab\\ ab \end{array}.$$

\begin{cons} La conjecture \ref{conj} est vraie pour $(\so_n, \so_{n_a}\times\so_{n_b})$ (BDI) avec $n_a\leqslant2$ ou $n_b\leqslant2$ ou $max(n_a,n_b)\leqslant 4$.
Dans le cas (CI), elle est vraie jusqu'au rang $7$; \emph{i.e.} pour $(\spn_{2n},\gl_{n})$ avec $n\leqslant 7$.
\end{cons}




 

\section{Une utilisation complète du lemme \ref{orbite}}
\label{pe1}
\subsection{\'Etude de $\g(e,1)$}
Nous avons jusqu'à présent écarté le cas où $\pp(e,0)$ contient des éléments nilpotents grâce à la proposition \ref{pd}. 
L'objet de cette section est d'utiliser le lemme \ref{orbite}, lorsque $\pp(e,1)$ est non nul dans les cas classiques qui nous restent. 
Notons que dans les cas BDI et CII, les presque $\pp$-distingués sont pairs. En particulier dans ces deux cas, pour un élément presque $\pp$-distingué $e$, 
on a $\pp(e,1)=\{0\}$. D'après le lemme \ref{norbite}, le lemme \ref{orbite} est impuissant pour éliminer la possibilité qu'un tel $e$ engendre une composante étrange. 
C'est pourquoi nous allons nous limiter aux cas AI et AII. 
Dans cette section on fixe donc $\g=\sln(V)$ et un élément nilpotent non nul $e$ que l'on inclut dans un $\Striplet$ $(e,h,f)$. 
Comme dans \cite{Ot2}, on note $(\lambda_i)_{i\in\N}$, la partition de $n$ correspondant au diagramme de Young de $e$. 
Alors, on sait qu'il existe une base \mbox{$\{e^a.v_i\mid a\leqslant\lambda_i-1; a,i\in \N, v_i\in V\}$} de $V$ (cf. \cite{TY} 19.2) pour laquelle $h$ est diagonal. 
Plus précisément, $$V_i:=\langle e^a.v_i\mid a\leqslant\lambda_i-1\rangle$$ est un $\sld$-module irréductible et en notant $V(k)$ 
le $h$-espace propre de $V$ associé à $k$, on a \begin{equation}\label{Vk} V(k)=\langle\{ e^a.v_i \mid 2a-\lambda_i+1=k\}\rangle.\end{equation}
 
Le lemme suivant nous indique des cas où il existe un élément $e_1\in \g^{e}$ tel que $e_1\notin\overline{\Od(e)}$.
\begin{lm}\label{lme1}
Supposons qu'il existe $i_1,i_2\in \N$ tels que $0\lnq\lambda_{i_1}=\lambda_{i_2}-1$. Alors, il existe un élément $e_1\in\bigoplus_{i\geqslant1}\g(e,i)$ tel que $e_1\notin\overline{\Od(e)}$.
\end{lm}
\begin{proof}
On note $e=\sum_{i\in \N} e_i$ où $e_i$ est défini par $\left\{\begin{array}{l}e \mbox{ sur } V_i\\ 0\mbox{ sur }\bigoplus_{j\neq i} V_j.\end{array}\right.$
Les éléments $e_i$ commutent deux à deux. Définissons maintenant l'élément 
\begin{equation} \label{11nov} e_{i_1,i_2}:\left\{\begin{array}{l} e^a.v_{i_1}\rightarrow e^{a+1}.v_{i_2} \mbox{ pour } 0\leqslant a\leqslant \lambda_{i_1}-1\\
e^a.v_{i_2}\rightarrow e^{a}.v_{i_1} \mbox{ pour } 0\leqslant a\leqslant \lambda_{i_1}-1.\\ e^a.v_j\rightarrow 0 \mbox{ sinon.}\end{array}\right.\end{equation}
On vérifie facilement que $e_{i_1}+e_{i_2}=e_{i_1,i_2}^2$, et donc que $e$ commute avec $e_1:=e_{i_1,i_2}+\sum_{i\neq i_1,i_2}e_i$. 
Par ailleurs le diagramme de $e_1$ est strictement supérieur à celui de $e$, donc $e_1\notin\overline{\Od(e)}$. 
Enfin, l'égalité \eqref{Vk} nous donne $e_{i_{1},i_{2}}.(V(k))\subseteq V(k+1)$ donc $e_1\in \g(e,1)\oplus\g(e,2)$.  
\end{proof}

D'après le lemme \ref{norbite}, on a $\bigoplus_{i\geqslant2}\g(e,i)\subset \overline{\Od(e)}$. 
Le lemme suivant peut alors être considéré comme une réciproque du lemme \ref{lme1}.
\begin{prop}\label{nlme1}
Supposons que pour tous $i,j\in \N$ tels que $\lambda_i,\lambda_j\neq 0$, on ait $\lambda_i-\lambda_j\neq 1$. Alors $\g(e,1)=\{0\}$.
\end{prop}
\begin{proof}
Nous allons montrer la contraposée en supposant l'existence de $x\in\g(e,1)\setminus\{0\}$.
Commençons par noter que d'après \eqref{Vk}, on a pour tous $j\in \N^*$, $k\in\N$:
\begin{equation}\label{Vkk}v_{j}\in V(-\lambda_{j}+1) \mbox{ et } \ker(e^{k})\cap V(-k+1)= \langle v_{i}\mid \lambda_{i}=k\rangle.
\end{equation}
Comme $x$ commute avec $e$, il existe un indice $i$ tel que $y=x.v_{i}\neq0$ 
et puisque $x\in\g(1,h)$, on a $y\in V(-\lambda_{i}+2)=e.V(-\lambda_{i})\oplus \ker (e^{\lambda_{i}-1})$. 
Soit $y_{1}+y_{2}$ la décomposition de $y$ dans cette somme directe. 
Soit $z\in V(-\lambda_{i})$ tel que $y_{1}=e.z$. On a $e^{\lambda_{i}+1}.z=x.(e^{\lambda_{i}}.v_{i})=0$ donc 
$z\in\langle v_{j}\mid \lambda_{j}=\lambda_{i}+1\rangle$. Si $z\neq0$, il existe nécessairement $j$ tel que $\lambda_{j}=\lambda_{i}+1$. 
Dans le cas contraire, $y_{1}=e.z=0$ et  
$0\neq y=y_{2}\in\ker(e^{\lambda_{i}-1})\cap V(-\lambda_{i}+2)=\langle v_{j}\mid \lambda_{j}=\lambda_{i}-1\rangle$. D'où l'existence de $j$ tel que $\lambda_{j}=\lambda_{i}-1$.
\end{proof}

Grâce au lemme \ref{norbite}, on voit que si $e$ est presque $\pp$-distingué et vérifie les hypothèses de la proposition précédente alors $\NN\cap\pp^e\subset\overline{\Od(e)}$. 
Le lemme \ref{orbite} ne peut donc rien nous apporter. 
Dans le cas contraire nous allons pouvoir utiliser ce lemme; c'est l'objet de la suite de cette section.

\subsection{Le cas AI}
On utilise les notations de la sous-section précédente. Soit $e$ presque $\pp$-distingué; les $\lambda_i$ sont donc deux à deux distincts. 
De plus, prenant en considération la proposition \ref{nlme1}, 
nous ferons l'hypothèse qu'il existe deux indices $i_1$ et $i_2$ tels que $0\lnq\lambda_{i_1}=\lambda_{i_2}-1$.
On sait qu'on peut alors supposer que la forme bilinéaire symétrique sur $V$ définissant la paire symétrique $(\g,\kk)$ est donnée par (cf. \cite[Lemme 1]{Ot2}) 
$$(e^a.v_i,e^b.v_j)=\left\{\begin{array}{l} 1\mbox{ si }i=j\mbox{ et }a+b+1=\lambda_i\\ 0 \mbox{ sinon}.\end{array}\right.$$
Il est facile de vérifier que pour tous $a,b,i,j,$ on a $(e_{1}.e^{a}.v_{i},e^b.v_{j})=(e^{a}.v_{i},e_{1}.e^b .v_{j})$ où
$e_1=e_{i_{1},i_{2}}+\sum_{j\neq i_{1},i_{2}}e_{i}$ est l'élément de la démonstration du lemme \ref{lme1}. Ceci implique que $e_{1}\in \pp$ et 
on a obtenu un élément de $\pp^e$ qui n'appartient pas à $\overline{\Od(e)}$. 
D'après le lemme \ref{orbite}, la sous-variété $\CC(e)$ n'est donc pas une composante irréductible de $\cnil$.

Concrètement, cela veut par exemple dire que l'élément correspondant au diagramme de Young $\Gamma_2$ de l'exemple \ref{exAI} ne peut pas engendrer de composante irréductible étrange. 
\begin{cons}La conjecture \ref{conj} est démontrée dans le cas AI de \mbox{rang 4}.\end{cons}
 
\subsection{Le cas AII}
Nous allons cette fois nous placer dans le cas AII. Soit $e$ presque $\pp$-distingué; donc si $\lambda_i\neq0$, 
il existe un unique indice $\beta(i)$ tel que $\lambda_{\beta(i)}=\lambda_i$. 
De plus, cf. proposition \ref{nlme1}, nous allons supposer qu'il existe $i_1$ et $i_2$ tels que 
$0\lnq\lambda_{i_1}=\lambda_{i_2}-1$.
On peut alors supposer que la forme bilinéaire antisymétrique sur $V$ définissant $(\g,\kk)$ est donnée par (cf. \cite[Lemme 1]{Ot2}) 
$$(e^a.v_i,e^b.v_j)=\left\{\begin{array}{l} \alpha(i)\mbox{ si }i=\beta(j)\mbox{ et }a+b+1=\lambda_i\\ 0 \mbox{ sinon}.\end{array}\right.$$ 
où $\alpha(i)\in\{\pm1\}$ et $\alpha(\beta(i))=-\alpha(i)$.
Quitte à permuter $i_1$ et $\beta(i_1)$, on peut supposer $\alpha(i_1)=\alpha(i_2)$. Avec des notations de la démonstration du lemme \ref{lme1} 
et en particulier \eqref{11nov}, on pose $e_{1}'=e_{i_{1},i_{2}}+e_{\beta(i_{1}),\beta(i_{2})}+\sum_{i\notin \{i_{1},i_{2},\beta(i_{1}), \beta(i_{2})\}} e_{i}$.
Comme dans le cas AI, on vérifie que c'est un élément de $\pp^e$ qui n'appartient pas à $\overline{\Od(e)}$. D'après le lemme \ref{orbite}, la sous-variété $\CC(e)$ n'est donc pas une composante irréductible de $\cnil$.

Ce qui précède implique par exemple que l'élément correspondant au diagramme de Young $\Gamma_4$ de l'exemple \ref{exAI} ne peut pas engendrer de composante irréductible étrange. 
\begin{cons}La conjecture \ref{conj} est démontrée dans le cas AII de \mbox{rang 3}.\end{cons}

\section{Les cas exceptionnels}
Dans les cas exceptionnels, nous disposons d'une classification des orbites d'éléments nilpotents, des centralisateurs des $\Striplet$s associés, et de leurs relations d'inclusion, cf. \cite{Do1}, \cite{Do2}, \cite{Do3}, \cite{Do4}. 
Nous allons donc pouvoir appliquer au cas exceptionnel les mêmes méthodes que précédemment visant à éliminer des possibilités, pour des orbites  nilpotentes, d'engendrer des composantes étranges. 
Finalement, nous démontrerons la conjecture dans tout les cas exceptionnels hormis EI.

\subsection{Utilisation de la proposition \ref{pd}}
Fixons une algèbre de Lie simple symétrique exceptionnelle $(\g,\kk)$ ainsi qu'une forme réelle du même type $\g_{\R}$.
Les tables de \cite{Do1} et \cite{Do2} nous donnent les orbites nilpotentes réelles dans $\g_{\R}$. En fait, la classification de D.Z.~Djokovic s'appuie sur une description des orbites complexes de $(\g,\kk)$ qu'il relie à la classification des orbites réelles par la correspondance de Kostant-Sekiguchi. 
Chaque ligne désigne donc une orbite complexe non nulle, que l'on repérera par son numéro donné dans la première colonne. Fixons un élément $e$ de cette orbite, et incluons-le dans un $\Striplet$ $(e,h,f)$ engendrant une algèbre de Lie $\sfr$. 
Le centralisateur réductif réel donné dans la dernière colonne de ces tables est en réalité déduit du centralisateur complexe $\g^{\sfr}$ (cf. \cite[\S15]{Do1}). 
Tous les calculs de D.Z.~Djokovic effectués dans le cas complexe restent vrais dans le cas d'un corps algébriquement clos de caractéristique zéro. 
Pour retrouver la paire réductive $(\g^{\sfr},\kk^{\sfr})$ complexe, il suffit d'appliquer la correspondance dans l'autre sens, sachant que $V_r$ désigne un tore de dimension $r$ dans $\pp$ et que $T_r$ désigne un tore de dimension $r$ dans $\kk$. 
On voit donc facilement si $\g^{\sfr}=\kk\oplus V_r$ ou non. Les tableaux suivants résultent de ces calculs dans les différents cas simples exceptionnels et donnent les orbites presque $\pp$-distingués. 
La première colonne donne le numéro de l'orbite tel qu'il apparaît dans \cite{Do1} et \cite{Do2}. La seconde colonne donne le type d'isomorphisme de la paire $(\g^{\sfr}, \kk^{\sfr})$. 
La troisième donne le défaut d'un élément de cette orbite.
\medskip

\begin{center}
\begin{tabular}{c c c}
\begin{tabular}{|c|c|c|}
\multicolumn{3}{c}{$E_{6(-26)}$ (cas EIV)}\\
\hline
1 & $(\so_7\oplus T_1,\so_7)$ & $1$\\
2 & $(G_2,G_2)$ & $0$\\ 
\hline
\end{tabular}
&
\begin{tabular}{|c|c|c|}
\multicolumn{3}{c}{ $F_{4(-20)}$ (cas FII)}\\
\hline
1 & $(\sln_4, \sln_4)$ & $0$\\
2 & $(G_2,G_2)$ & $0$\\
\hline
\end{tabular}
&
\begin{tabular}{|c|c|c|}
\multicolumn{3}{c}{$G_{2(2)}$ (cas GI)}\\
\hline
3 & $(0,0)$ & $0$\\
4 & $(0,0)$ & $0$\\
5 & $(0,0)$ & $0$\\
\hline
\end{tabular}
\end{tabular}

\medskip

\begin{tabular}{c c c}
\begin{tabular}{|c|c|c|}
\multicolumn{3}{c}{$E_{6(-14)}$ (cas EIII)}\\
\hline
3 & $(\so_7\oplus T_1,\so_7\oplus T_1)$ & $0$\\
4 & $(\so_7\oplus T_1,\so_7\oplus T_1)$ & $0$\\ 
7 & $(\sln_3\oplus T_1,\sln_3\oplus T_1)$ & $0$\\
8 & $(\sln_3\oplus T_1,\sln_3\oplus T_1)$ & $0$\\
9 & $(G_2,G_2)$ & $0$\\
10 & $(\so_5\oplus T_1,\so_5\oplus T_1)$ & $0$\\
11 & $(\so_5\oplus T_1,\so_5\oplus T_1)$ & $0$\\
12 & $(\sln_2\oplus T_1,\sln_5\oplus T_1)$ & $0$\\
\hline
\end{tabular}
&
\begin{tabular}{|c|c|c|}
\multicolumn{3}{c}{$E_{6(6)}$ (cas EI)}\\
\hline
12 & $(T_2,T_1)$ & $1$\\
16 & $(T_1,0)$ & $1$\\
17 & $(T_1,0)$ & $1$\\
18 & $(0,0)$ & $0$\\
19 & $(0,0)$ & $0$\\
20 & $(0,0)$ & $0$\\
21 & $(T_1,0)$ & $1$\\
22 & $(0,0)$ & $0$\\
23 & $(T_2,0)$ & $2$\\
\hline
\end{tabular}
&
\begin{tabular}{|c|c|c|}
\multicolumn{3}{c}{$F_{4(4)}$ (cas FI)}\\
\hline
6 & $(\sln_3,\sln_3)$ & $0$\\
16 & $(0,0)$ & $0$\\
17 & $(0,0)$ & $0$\\
18 & $(0,0)$ & $0$\\
19 & $(\sln_2,\sln_2)$ & $0$\\
22 & $(0,0)$ & $0$\\
23 & $(0,0)$ & $0$\\
24 & $(0,0)$ & $0$\\
25 & $(0,0)$ & $0$\\
26 & $(0,0)$ & $0$\\
\hline
\end{tabular}
\end{tabular}

\begin{tabular}{c c}
\begin{tabular}{|c|c|c|}
\multicolumn{3}{c}{$E_{7(-25)}$ (cas EVII)}\\
\hline
6 & $(F_4,F_4)$ & $0$\\
7 & $(F_4,F_4)$ & $0$\\
11 & $(\sln_4\oplus T_1,\sln_4\oplus T_1,)$ & $0$\\
12 & $(\sln_4\oplus T_1,\sln_4\oplus T_1,)$ & $0$\\
16 & $(\so_7,\so_7)$ & $0$\\
17 & $(\so_7,\so_7)$ & $0$\\
18 & $(\so_7,\so_7)$ & $0$\\
19 & $(\so_7,\so_7)$ & $0$\\
20 & $(\sln_3\oplus T_1,\sln_3\oplus T_1)$ & $0$\\
21 & $(G_2,G_2)$ & $0$\\
22 & $(G_2,G_2)$ & $0$\\
\hline
\end{tabular}
&
\begin{tabular}{|c|c|c|}
\multicolumn{3}{c}{$E_{8(-24)}$ (cas EIX)}\\
\hline
6 & $(E_6,E_6)$ & $0$\\
18 & $(\so_8,\so_8)$ & $0$\\
19 & $(\so_8,\so_8)$ & $0$\\
21 & $(F_4,F_4)$ & $0$\\
23 & $(\so_5\oplus T_1,\so_5\oplus T_1)$ & $0$\\
24 & $(\sln_5, \sln_5)$ & $0$\\
26 & $(\sln_3\oplus T_1,\sln_3\oplus T_1)$ & $0$\\
27 & $(\sln_4,\sln_4)$ & $0$\\
28 & $(2\sln_2,2\sln_2)$ & $0$\\
30 & $(G_2,G_2)$ & $0$\\
31 & $(G_2,G_2)$ & $0$\\
32 & $(\so_7,\so_7)$ & $0$\\
33 & $(\so_7,\so_7)$ & $0$\\
34 & $(2\sln_2,2\sln_2)$ & $0$\\
35 & $(\sln_3,\sln_3)$ & $0$\\
36 & $(G_2,G_2)$ & $0$\\
\hline
\end{tabular}
\end{tabular}

\begin{tabular}{c c}
\begin{tabular}{|c|c|c|}
\multicolumn{3}{c}{$E_{7(7)}$ (cas EV)}\\
\hline
16 & $(G_2,G_2)$ & $0$\\
17 & $(G_2,G_2)$ & $0$\\
39 & $(\sln_2,\sln_2)$ & $0$\\
40 & $(\sln_2,\sln_2)$ & $0$\\
48 & $(T_2,T_2)$ & $0$\\
49 & $(T_2,T_2)$ & $0$\\
50 & $(T_2,0)$ & $2$\\
55 & $(\sln_2,\sln_2)$ & $0$\\
56 & $(\sln_2,\sln_2)$ & $0$\\
67 & $(0,0)$ & $0$\\
68 & $(0,0)$ & $0$\\
69 & $(0,0)$ & $0$\\
70 & $(0,0)$ & $0$\\
76 & $(0,0)$ & $0$\\
77 & $(0,0)$ & $0$\\
78 & $(0,0)$ & $0$\\
79 & $(0,0)$ & $0$\\
80 & $(T_1,T_1)$ & $0$\\
81 & $(T_1,0)$ & $1$\\
85 & $(0,0)$ & $0$\\
86 & $(0,0)$ & $0$\\
87 & $(0,0)$ & $0$\\
88 & $(0,0)$ & $0$\\
89 & $(0,0)$ & $0$\\
90 & $(0,0)$ & $0$\\
91 & $(0,0)$ & $0$\\
92 & $(0,0)$ & $0$\\
93 & $(0,0)$ & $0$\\
94 & $(0,0)$ & $0$\\
\hline
\end{tabular}
&
\begin{tabular}{|c|c|c|}
\multicolumn{3}{c}{$E_{8(8)}$ (cas EVIII)}\\
\hline
14 & $(G_2,G_2)$ & $0$\\
15 & $(G_2,G_2)$ & $0$\\
34 & $(\sln_3,\sln_3)$ & $0$\\
42 & $(\sln_2\oplus T_1,\sln_2\oplus T_1)$ & $0$\\
45 & $(2\sln_2,2\sln_2)$ & $0$\\
51 & $(\sln_3,\sln_3)$ & $0$\\
67 & $(0,0)$ & $0$\\
68 & $(0,0)$ & $0$\\
69 & $(0,0)$ & $0$\\
70 & $(2\sln_2,2\sln_2)$ & $0$\\
79 & $(T_1,T_1)$ & $0$\\
80 & $(T_1,T_1)$ & $0$\\
81 & $(T_1,0)$ & $1$\\
84 & $(T_1,T_1)$ & $0$\\
85 & $(T_1,0)$ & $1$\\
87 & $(T_1,T_1)$ & $0$\\
88 & $(T_1,0)$ & $1$\\
91 & $(0,0)$ & $0$\\
92 & $(0,0)$ & $0$\\
93 & $(T_1,T_1)$ & $0$\\
94 & $(T_1,T_1)$ & $0$\\
95 & $(T_1,0)$ & $1$\\
98 & $(0,0)$ & $0$\\
99 & $(0,0)$ & $0$\\
101 & $(0,0)$ & $0$\\
102 & $(0,0)$ & $0$\\
104 & $(0,0)$ & $0$\\
105 & $(0,0)$ & $0$\\
106 & $(0,0)$ & $0$\\
107 & $(0,0)$ & $0$\\
109 & $(0,0)$ & $0$\\
110 & $(0,0)$ & $0$\\
111 & $(0,0)$ & $0$\\
112 & $(0,0)$ & $0$\\
113 & $(0,0)$ & $0$\\
114 & $(0,0)$ & $0$\\
115 & $(0,0)$ & $0$\\
\hline
\end{tabular}
\end{tabular}

\begin{tabular}{c c}
\begin{tabular}{|c|c|c|}
\multicolumn{3}{c}{$E_{6(2)}$ (cas EII)}\\
\hline
6& $(2\sln_2,2\sln_2)$ & $0$\\
12 & $(\sln_2\oplus T_1,\sln_2\oplus T_1)$ & $0$\\
13 & $(\sln_2\oplus T_1,\sln_2\oplus T_1)$ & $0$\\
20 & $(T_2,T_2)$ & $0$\\
21 & $(T_2,T_2)$ & $0$\\
22 & $(T_2,T_1)$ & $1$\\
23 & $(\sln_3, \sln_3)$ & $0$\\
25 & $(\sln_2\oplus T_1, \sln_2\oplus T_1)$ & $0$\\
27 & $(T_1,T_1)$ & $0$\\
28 & $(T_1,T_1)$ & $0$\\
29 & $(T_1,T_1)$ & $0$\\
30 & $(T_1,T_1)$ & $0$\\
32 & $(0,0)$ & $0$\\
33 & $(0,0)$ & $0$\\
34 & $(T_1,T_1)$ & $0$\\
35 & $(T_1,T_1)$ & $0$\\
36 & $(0,0)$ & $0$\\
37 & $(0,0)$ & $0$\\
\hline
\end{tabular}
&
\begin{tabular}{|c|c|c|}
\multicolumn{3}{c}{$E_{7(-5)}$ (cas EVI)}\\
\hline
6 & $(\sln_6,\sln_6)$ & $0$\\
14 & $(G_2+\sln_2,G_2+\sln_2)$ & $0$\\
19 & $(3\sln_2,3\sln_2)$ & $0$\\
20 & $(3\sln_2,3\sln_2)$ & $0$\\
22 & $(\spn_6,\spn_6)$ & $0$\\
24 & $(\sln_2\oplus T_1,\sln_2\oplus T_1)$ & $0$\\
25 & $(\sln_3\oplus T_1,\sln_3\oplus T_1)$ & $0$\\
27 & $(T_2,T_2)$ & $0$\\
28 & $(\sln_2\oplus T_1,\sln_2\oplus T_1)$ & $0$\\
29 & $(\sln_2,\sln_2)$ & $0$\\
31 & $(\sln_2,\sln_2)$ & $0$\\
32 & $(\sln_2,\sln_2)$ & $0$\\ 
33 & $(2\sln_2,2\sln_2)$ & $0$\\
34 & $(2\sln_2,2\sln_2)$ & $0$\\
35 & $(\sln_2,\sln_2)$ & $0$\\
36 & $(T_1,T_1)$ & $0$\\
37 & $(\sln_2,\sln_2)$ & $0$\\
\hline
\end{tabular}
\end{tabular}
\end{center}

Notons que l'orbite 15 de EVIII n'apparaît pas dans \cite{PT}, elle est néanmoins bien $\pp$-distinguée. 
Notons aussi que l'orbite 1 de EIV n'est pas $\pp$-distinguée suite à une erreur dans la table VII de \cite{Do2} mentionnée par King dans \cite{Ki}.
\begin{cons}Les algèbres de Lie simples symétriques GI, FI, FII, EIII, EVI, EVII et EIX ne possèdent pas d'élément presque $\pp$-distingué non $\pp$-distingué 
donc en vertu de la proposition \ref{pd}, la conjecture \ref{conj} est démontrée dans ces cas-là. \end{cons}

\subsection{Quelques réductions}
Dans une algèbre de Lie simple symétrique donnée, on note $\Od_i$ l'orbite de numéro $i$ des tables de D.Z.~Djokovic. 
\'Enumérons quelques réductions d'éléments presque $\pp$-distingués.
\nopagebreak[4]
\vspace{-0.15cm}
\begin{center}
\begin{tabular}{|c|c|c|c|c|}
\hline
Reference&Type de $\g$&  n° & $\dim \Od$ & $\delta(e)$\\
\hline
\multirow{2}{2cm}{\cite{Do3}} &\multirow{2}{2cm}{EII} & 22 & 29 & 1\\
\cline{3-5}
& &24 & 30 & 0\\
\hline
\hline
\multirow{4}{2cm}{\cite{Do4}} &\multirow{4}{2cm}{EV} &50 & 52 & 2\\
\cline{3-5}
& &54 & 53 & 1\\
\cline{3-5}
\cline{3-5}
& &81 & 59 & 1\\
\cline{3-5}
& &85 & 60 & 0\\
\hline
\hline
\multirow{6}{2cm}{\cite{Do5} Table 2 et Fig.2} &\multirow{6}{2cm}{EVIII} & 81 & 107 & 1\\
\cline{3-5}
& &84 & 108 & 0\\
\cline{3-5}
\cline{3-5}
& &88 & 109 & 1\\
\cline{3-5}
& &91 & 110 & 0\\
\cline{3-5}
\cline{3-5}
& &95 & 111 & 1\\
\cline{3-5}
& &98 & 112 & 0\\
\hline
\hline
\multirow{4}{2cm}{\cite{Do3} Fig.2} &\multirow{4}{2cm}{EI} &21 & 34 &1\\
\cline{3-5}
& &18 & 35 &0\\
\cline{3-5}
\cline{3-5}
& &17 & 32 & 1\\
\cline{3-5}
& &22 & 33 & 0\\
\hline
\end{tabular}
\end{center}

Il ne reste que 5 orbites sur lesquelles nous ne pouvons pas nous prononcer pour l'instant; elles ne possèdent pas de réduction. Ce sont les suivantes.

\begin{center}
\begin{tabular}{|c|c|c|c|}
\hline
Type & n° de $\Od$ & \begin{tabular}{c} $K$-diagramme de Dynkin \\ d'une caractéristique \\de $\Od$\end{tabular} & \begin{tabular}{c}$G$-diagramme de Dynkin \\ d'une caractéristique \\de $\Od$ \end{tabular} \\
\hline
EIV & 1& 0001 & 100001 \\
EI & 16 & 1111& 111011 \\
EVIII & 85 & 11111111 & 10010101\\
EI & 12 & 2002 & 000200  \\
EI & 23 & 0020 & 000200  \\
\hline
\end{tabular}
\end{center}

\subsection{Utilisation du lemme \ref{orbite}}
\label{pe1bis}
Nous allons montrer que les trois premières orbites du tableau précédent n'engendrent pas de composante étrange, grâce au lemme \ref{orbite}. 
Comme dans la section \ref{pe1}, on étudie le $h$-espace propre associé à la valeur propre $1$.

Soit $\g= E_{6(-26)}$ (cas EIV). On s'intéresse à l'orbite $\Od_1$ (cf. \cite{Do2}). 
On considère un $\Striplet$ normal standard $(e,h,f)$ avec $e\in \Od_1$. 
On peut calculer facilement les dimensions de $\kk(i,h)$ et $\pp(i,h)$ pour les différents entiers $i$ et on obtient (cf. \cite{JN}):

\begin{center}
\begin{tabular}{c||c|c}
$i$ & $\dim \kk(i,h)$ & $\dim \pp(i,h)$\\
\hline
\hline
$\gnq 2$ & 0 & 0\\
\hline
2 & 7& 1\\
\hline
1 & 8 & 8\\
\hline
0 & 22 & 8\\
\end{tabular}  
\end{center}

On en déduit que tout élément non nul $e$ de $\pp(2,h)$ est régulier sous l'action de $K(0,h)$ dans $\pp(2,h)$ et appartient donc à $\Od_{1}$ (cf. \cite[Lemme 2.2.9]{Ka}). 
De la même façon, comme $h'=2h$ est une caractéristique pour l'orbite régulière $\Od_{2}$, tout élément régulier sous l'action de $K(0,h')=K(0,h)$ dans $\pp(2,h')=\pp(1,h)$ appartient à $\Od_2$. 
Mais comme $\kk(3,h)=\{0\}$, on a $\pp(1,h)=\pp(e,1)$. 
Soit donc $e'$ régulier dans $\pp(e,1)$; on a $e'\in \Od_2$ et $[e',e]=0$, d'où $\mathcal N\cap\pp^e \not\subset \overline{\Od_1}$, et par le lemme \ref{orbite}, $e$ ne peut pas engendrer de composante étrange.
\newline

\begin{prop} Si $(\g,\kk)$ est de type EI, il existe $e\in\Od_{16}$ et $e_{1}\in\Od_{18}$ tels que $[e, e_{1}]=0$.\item 
Si $(\g,\kk)$ est de type EVIII, il existe $e\in\Od_{85}$ et $e_{1}\in\Od_{109}$ tels que $[e, e_{1}]=0$.
\end{prop}
\begin{proof}
Soit $\g$ de type $E_{6}$ (resp. $E_{8}$). Nous allons fixer une sous algèbre de Cartan $\h$ de $\g$, une base de Chevalley $(X_{\alpha_{i}})_{\alpha_{i}\in R(\g,\h)}$ 
compatible avec le système de racine $R(\g,\h)$, et une involution explicite de type EI (resp. EVIII) telle que $\h\cap\kk$ soit une sous-algèbre de Cartan de $\kk$.
Pour simplifier les notations, on notera parfois $X_{i}$ à la place de $X_{\alpha_{i}}$. 
Soit une caractéristique $h\in\h\cap\kk$ de l'orbite $\Od_{16}$ (resp. $\Od_{85}$). Il se trouve que $h'=2h$ est alors une caractéristique de l'orbite $\Od_{18}$
(resp. $\Od_{109}$). Avec les notations de la section~1.2, on rappelle que l'on a $\kk(0,h)=\kk(0,h')$ et $\pp(2,h')=\pp(1,h)$. 
De plus, d'après \cite[Lemme~2.2.9]{Ka}, $\Od_{16}\cap\pp(2,h)$ (resp. $\Od_{85}\cap\pp(2,h)$) est exactement la $K(0,h)$-orbite régulière dense de $\pp(2,h)$, 
tandis que $\Od_{18}\cap\pp(1,h)$ (resp. $\Od_{109}\cap\pp(1,h)$) est la $K(0,h)$-orbite régulière dense de $\pp(1,h)$. 
Grâce aux tables de \cite{JN}, on voit que pour tout $e\in\Od_{16}\cap\pp(2,h)$, on a $\dim \pp(e,1)=1$. 
Nous choisirons donc un élément régulier $e\in\pp(e,2)$ et montrer qu'il commute avec un élément $e_{1}$, régulier dans $\pp(1,h)$.

Commençons par le premier cas. On fixe une base $(\alpha_{i})_{1\leqslant i\leqslant 6}$ telle que le diagramme de Dynkin de $\g$ s'écrive
$$\begin{array}{rrrrr} \alpha_{1}&\alpha_{3}&\alpha_{4}&\alpha_{5}&\alpha_{6} \\ &&\alpha_{2}&&\end{array}.$$
On choisit la numérotation des racines donnée par \cite[Table 3]{Do6} ainsi que la base de Chevalley de \cite[\S4]{Do6} caractérisée par \cite[Table 13]{Do6}.
Ainsi, on a par exemple $\alpha_{8}=\alpha_{2}+\alpha_{4}$ et $[X_{2}, X_{4}]=-X_{8}$. On fixe l'involution $\tau$, définie dans \cite{Do2} 
(dans \cite{Do2} les racines $\alpha_{2},\alpha_{3},\alpha_{4}$ sont notées respectivement $\alpha_{4},\alpha_{2},\alpha_{3}$), 
permutant les racines $\alpha_1$ et $\alpha_6$, $\alpha_3$ et $\alpha_5$ et laissant invariants $\alpha_2$ et $\alpha_4$. 
Pour trouver une involution $\theta$ de type EI, on peut poser (cf. \cite{Do2})
$$\theta\big(X_{\sum_{i=1}^6 k_i \alpha_i}\big)=(-1)^{k_2}X_{\sum_{i=1}^6 k_i \tau(\alpha_i)}.$$
Toujours selon \cite{Do2} on définit $$\begin{array}{c} \beta_{1}=(\alpha_{1})_{\mid\h\cap\kk}=(\alpha_{6})_{\mid\h\cap\kk}, \; \;
\beta_{2}=(\alpha_{3})_{\mid\h\cap\kk}=(\alpha_{5})_{\mid\h\cap\kk}, \; \;\beta_{3}=(\alpha_{4})_{\mid\h\cap\kk},\\
\beta_{4}=(\alpha_{2})_{\mid\h\cap\kk}, \qquad \beta_{0}=-2\beta_{1}-3\beta_{2}-2\beta_{3}-\beta_{4}.\end{array}$$ 
De \cite[Table VIII]{Do2}, on déduit que l'élément $h$ vérifiant $\beta_1(h)=\beta_2(h)=\beta_3(h)=\beta_0(h)=1$ est une caractéristique de $\Od_{16}$.
Si on écrit le $G$-diagramme de Dynkin de $h$, on trouve $$\begin{array}{rrrrr} 1&1&\!\!\!\!1&1&1 \\ &&\!\!\!\!-8&&\end{array}.$$
Ceci nous permet d'obtenir les bases suivantes pour $\kk(0,h)$, $\pp(1,h)$, $\pp(2,h)$ et $\kk(3,h)$:
\begin{center}
\begin{tabular}{|c|c|}
\hline
$\kk(0,h)$ & $H_{1}+H_{6}; H_{3}+H_{5}; H_{2}; H_{4}$\\
$\pp(1,h)$ & $A_1:=X_1-X_6;A_2:=X_3-X_5; A_3:=X_{-32}+X_{-33}; A_4:=X_{35}$\\
$\pp(2,h)$ & $B_1:=X_7-X_{11}; B_2:=X_9-X_{10}; B_3:=X_{-29}+X_{-31}; B_4:=X_{-30}$\\
$\kk(3,h)$ & $C_1:=X_{12}+X_{16}; C_2:=X_{-26}-X_{-28};C_3:=X_{15}$\\
\hline
\end{tabular}
\end{center}
Comme $\kk(0,h)\subset \h$, un élément $x\in\h$ agit par multiplication par un scalaire sur chaque $A_{i}$ et chaque $B_{j}$. 
En particulier, un élément $e=\sum b_{i}B_{i}$ (resp. $e_{1}=\sum a_{i}A_{i}$) est régulier dans $\pp(2,h)$ (resp. $\pp(1,h)$) si et seulement si 
$[\kk(0,h),e]=\pp(2,h)$ (resp. $[\kk(0,h), e_{1}]=\pp(1,h)$) ce qui équivaut au fait que
les coefficients $b_{i}$ (resp. $a_{i}$) sont tous non nuls. On choisit $e=B_{1}+B_{2}+B_{3}+B_{4}$.
Ecrivons la table de multiplication dans la matrice ci-dessous, où l'élément d'indice $(i,j)$ correspond à $[A_i,B_j]$:
$$\big([A_i,B_j]\big)_{i,j}=\left(\begin{array}{c c c c} 0 &C_1& 0& C_2\\0 &-2C_3 &-C_2 &0\\ C_2 & 0& 0 &0\\ 0 & 0 & C_1 & C_3\end{array}\right)$$
D'après ce qui précède, l'élément $e_{1}=-2A_{1}+A_{2}+3A_{3}+2A_{4}$ est régulier dans $\pp(1,h)$, appartient à $\Od_{18}$ et commute avec $e$. 

Le second cas où $(\g,\kk)$ est de type EVIII se traite de façon similaire. On fixe la base de $R(\g,\h)$:
$$\begin{array}{rrrrrrr} \alpha_{1}&\alpha_{3}&\alpha_{4}&\alpha_{5}&\alpha_{6}& \alpha_{7} &\alpha_{8} \\ &&\alpha_{2}&& & &\end{array}.$$
La numérotation des racines et la base de Chevalley sont données par \cite[Table 1, Table 6]{Do7}. Conformément à \cite{Do1}, on pose
$$\theta\big(X_{\sum_{i=1}^8 k_i \alpha_i}\big)=(-1)^{k_1}X_{\sum_{i=1}^8 k_i \alpha_i}.$$
On peut choisir la caractéristique $h$ de $\Od_{85}$ de sorte que son $G$-diagramme de Dynkin soit
$$\begin{array}{rrrrrrr} -14&1&1&1&1&1&1 \\ &&1&&\end{array}.$$
On obtient facilement des bases pour les espaces de poids suivant:
\begin{center}
\begin{tabular}{|c|c|}
\hline
$\kk(0,h)$ & $H_{1}; H_{2}; H_{3}; H_{4}; H_{5}; H_{6}; H_{7}; H_{8}$\\
\hline
$\pp(1,h)$ & $A_1:=X_{93};A_2:=X_{94}; A_3:=X_{95}; A_4:=X_{96}$\\
& $A_5:=X_{-85};A_6:=X_{-86}; A_7:=X_{-87}; A_8:=X_{-88}$\\
\hline
$\pp(2,h)$ & $B_1:=X_{98}; B_2:=X_{99}; B_3:=X_{100}; B_4:=X_{-80}$\\
&$B_5:=X_{-81}; B_6:=X_{-82}; B_7:=X_{-83}; B_8:=X_{-84}$\\
\hline
$\kk(3,h)$ & $C_1:=X_{17}; C_2:=X_{18}; C_3:=X_{19}; C_4:=X_{20}$\\
&$C_5:=X_{21}; C_6:=X_{22}; C_7:=X_{-118}$\\
\hline
\end{tabular}
\end{center}
Ici encore $\kk(0,h)\subseteq\h$, donc les orbites régulières de $\pp(1,h)$ et $\pp(2,h)$ se décrivent de la même façon que dans le cas EI.
La table de multiplication: $\pp(1,h)\times\pp(2,h)\rightarrow\kk(3,h)$ est la suivante:
$$\big([A_i,B_j]\big)_{i,j}=\left(\begin{array}{c c c c c c c c} 0&0&0& -C_3& 0 & -C_1 &0&0\\ 0&0&0&0&C_4 &0 &C_2&0\\
0&0&0&0&C_5&0&0&C_2\\ 0&0&0&0&0&C_6&C_5&C_4\\ 0&0 &-C_6 & 0&0&0&-C_7&0\\ C_3 & 0&-C_5&0&0&-C_7&0&0\\ C_1&0&0&-C_7&0&0&0&0\\ 0&-C_2&0&0&0&0&0&0\end{array}\right)$$
Finalement, on peut poser $$e=\sum_{i=1}^8 B_{i}\mbox{ et }e_{1}=A_{1}+2A_{2}+3A_{3}-2A_{4}-2A_{5}+A_{6}+A_{7}+5A_{8}.$$ 
Ces deux éléments commutent et appartiennent respectivement à $\Od_{85}$ et $\Od_{109}$.\end{proof}
 
\begin{cor}
Si $(\g,\kk)$ est de type EI, $\CC(\Od_{18})$ n'engendre pas de composante étrange. Même conclusion pour $\Od_{85}$ si $(\g,\kk)$ est de type EVIII.
\end{cor}
\begin{proof}
Il suffit d'appliquer la proposition précédente et le lemme~\ref{orbite}.\end{proof}

Terminons cette section en donnant le nombre de composantes irréductibles de $\cnil$ dans les différents cas exceptionnels.

\begin{center}
\begin{tabular}{|c|c|c|c|c|c|c|c|c|c|c|c|}
\hline
GI & FI & FII & EI & EII & EIII & EIV & EV & EVI & EVII & EVIII & EIX\\
\hline
3& 10 & 2 & 4 à 6 & 17 & 8 & 1 & 27 & 17 & 11 &33 & 16\\
\hline
\end{tabular}
\end{center}

\section*{Conclusion}
La conjecture \ref{conjec} est démontrée dans les cas suivants:
AIII, CII, DIII, EII, EIII, EIV, EV, EVI, EVII, EVIII, EIX, FI, FII, GI
et dans les cas
\begin{itemize}
\item[AI] $(\sln_n, \so_n)$ avec $n\leqslant 5$ (\emph{i.e.} en rang $\leqslant 4$)
\item[AII] $(\sln_{2n}, \spn_{2n})$ avec $n\leqslant3$ (\emph{i.e} en rang $\leqslant 3$).
\item[BDI] $(\so_n, \so_p\times\so_q)$ avec $p\leqslant2$ ou $q\leqslant2$ ou $max(p,q)\leqslant 4$. (En particulier, elle est vraie en rang $\leqslant2$)
\item[CI] $(\spn_{2n},\gl_n)$ avec $n\leqslant7$. (\emph{i.e.} en rang $\leqslant7$)
\end{itemize}
De plus, les méthodes de réduction et d'étude de $\pp(e,1)$ introduites dans les sections \ref{reduction} et \ref{pe1} montrent qu'un certain nombre d'éléments presque $\pp$-distingués ne fournissent pas de composante étrange.

Nous avons exploité au maximum les deux outils dont nous nous sommes dotés: d'une part le lemme \ref{orbite} et d'autre part la réduction dont le principe est contenu dans le lemme \ref{xx}.
En effet, les seules orbites pour lesquelles nous ne savons rien dire sont des orbites $\Od$ contenant des éléments $e$ presque $\pp$-distingués, non $\pp$-distingués, vérifiant $\pp^e\cap \NN\subseteq \overline{\Od}$, et n'ayant pas de réduction.

\section{Appendice: Orbites $\pp$-self-large}
Le lemme \ref{orbite} nous a invité à déterminer les éléments nilpotents $e$ vérifiant $\NN(\pp^{e})\subset \overline{K.e}$. 
Dans le cas des algèbres de Lie, D.~Panyushev a dernièrement nommé ces éléments \emph{self-large} dans \cite{Pa4}. 
Ils sont caractérisés de la façon suivante: ce sont les éléments presque distingués $e$ vérifiant $\g(e,1)=0$. 
N'ayant découvert ces travaux que récemment, nous allons pré\-sen\-ter dans cet appendice quelques résultats inspirés de \cite{Pa4} qui permettent de dé\-ter\-mi\-ner 
des éléments similaires dans les algèbres de Lie symétriques. 

\begin{defi}
Un élément nilpotent $e\in\pp$ vérifiant $\NN(\pp^{e})\subset \overline{K.e}$ est dit \emph{$\pp$-self-large}. 
L'orbite $K.e$ est alors également appelée $\pp$-self-large  
\end{defi}

On a montré dans la section~1.4, que seuls des éléments $\pp$-self-large pouvaient engendrer une composante irréductible de $\cnil$.
On a en particulier établi les implications suivantes:
$$\pp\textrm{-distingué}\Rightarrow\pp\textrm{-self-large}\Rightarrow \textrm{presque }\pp\textrm{-distingué}.$$

Comme précédemment, on fixe un $\Striplet$ $(e,h,f)$, ce qui nous donne les graduations suivantes de $\wfr=\g, \kk$ ou $\pp$:
$$\wfr=\bigoplus_{i\in\Z}\wfr(i,h); \qquad \wfr^{e}=\bigoplus_{i\geqslant 0} \wfr(e,i); \qquad \wfr^{f}=\bigoplus_{i\leqslant 0} \wfr(f,i).$$ 
On sait que $\pp(e,0)$ agit (via l'action adjointe de $\g$) sur $\g(f,-1)$. Si $e$ est presque $\pp$-distingué, $\pp(e,0)$ est un tore, 
on peut alors considérer $\mathfrak{X}(\pp(e,0))$ l'ensemble des poids non-nuls du $\pp(e,0)$-module $\g(f,-1)$. 
On notera la décomposition en espaces de poids comme ceci:
$$\g(f,-1)=\bigoplus_{\gamma\in\mathfrak{X}(\pp(e,0))}V_{\gamma}.$$

Le but de la proposition \ref{sl} est de donner un analogue faible de \cite[Théorème 2.1]{Pa4} pour les algèbres de Lie symétriques.
Il se trouve que, combinée aux résultats des sections précédentes, cette proposition va s'avérer suffisante pour décrire toutes les orbites $\pp$-self-large des algèbres de Lie symétriques.
Nous avons tout d'abord besoin de deux lemmes. On rappelle que $L$ désigne la forme de Killing sur $\g$.

\begin{lm}[D.~Panyushev \cite{Pa4}]\label{mmmcor}
L'application $$\Phi: \left\{\begin{array}{rcl}\g(f,-1)\times\g(f,-1)&\rightarrow&\K\\
(\xi,\eta)&\mapsto & L(e,[\xi,\eta])\end{array}\right.$$
 est une forme bilinéraire antisymétrique non dégénérée $\g(e,0)$-invariante.
\end{lm}
\smallskip

\begin{lm}
L'automorphisme $\theta$ induit une bijection entre $V_{\gamma}$ et $V_{-\gamma}$. 
L'application $\Phit: (\xi,\eta)\mapsto\Phi(\xi,\theta(\eta))$ est une forme bilinéaire symétrique non dégénérée. 
Elle reste non dégénérée sur chaque sous-espace $V_{\gamma}$.
\end{lm}
\begin{proof}
Comme $h\in\kk$ et $f\in\pp$, l'automorphisme $\theta$ induit une bijection de $\g(f,-1)$. Soit $\xi$ un élément de $V_{\gamma}$. 
Pour tout $t\in\pp(e,0)$, on a
$$\label{dpt}[t,\theta(\xi)]=\theta([\theta(t),\xi])=\theta(-\gamma(t)\xi)=-\gamma(t)\theta(\xi).$$
Ceci prouve la première affirmation.

Comme $\theta$ induit une bijection de $\g(f,-1)$, le fait que $\Phit$ soit non-dégénérée est une conséquence du lemme \ref{mmmcor}.
Vérifions que $\Phit$ est symétrique:
\begin{eqnarray*}
L(e,[\xi,\theta(\eta)])&=&L(\theta(e),\theta([\xi,\theta(\eta)]))\\
&=&L(-e,[\theta(\xi),\eta])\\
&=&L(e,[\eta,\theta(\xi)]).
\end{eqnarray*}
Prouvons enfin la dernière affirmation du lemme. Soit $\xi\in V_{\gamma}$ et $\eta\in V_{\mu}$, alors $\theta(\eta)\in V_{-\mu}$ et pour tout $t\in\pp(e,0)$, on a
 $$(\gamma(t)-\mu(t))\Phit(\xi,\eta)=\Phi([t,\xi],\theta(\eta))+\Phi(\xi,[t,\theta(\eta)])=0.$$
donc $\Phit(\xi,\eta)=0$ si $\gamma\neq\mu$.
Ceci montre que $\Phit$ est non-dégénérée sur $V_{\gamma}$. 
\end{proof}

\begin{prop}\label{sl}
Si $e$ est presque $\pp$-distingué tel que $\g(f,-1)^{\pp(e,0)}=\{0\}$ et $\pp(e,1)\neq\{0\}$, alors $e$ n'est pas $\pp$-self-large.
\end{prop}
\begin{proof}
Notons tout d'abord que, d'après les hypothèses, $0\notin\mathfrak{X}(\pp(e,0))$.
Fixons $\mu\in\mathfrak{X}(\pp(e,0))$ et soit $\xi\in V_{\mu}$ tel que $\Phit(\xi,\xi)\neq0$. Puisque $\mu\neq 0$, il existe $t\in\pp(e,0)$ tel que $[t,\xi]=\xi$ et $[t,\theta(\xi)]=-\theta(\xi)$.
Un calcul facile donne 
$$L([[e,(\xi+\theta(\xi))],(\xi+\theta(\xi))],t)=2L(e,[\xi,\theta(\xi)])\neq 0,$$
ce qui montre en particulier que $[[e,\xi+\theta(\xi)],\xi+\theta(\xi)]\neq0$. Donc $z=\xi+\theta(\xi)$ est un élément de $\kk(f,-1)$ satisfaisant $[[e,z],z]\neq 0$. 
Finalement, par \cite[Lemme 2.3]{Pa4}, $G.(e+[z,e])$ est une orbite strictement plus grande que $G.e$, ce qui implique notamment que $K.(e+[z,e])\not\subseteq\overline{K.e}$.
Pour conclure, il est facile de vérifier que $e+[z,e]\in\pp^{e}$. 
\end{proof}

\begin{Rq}\label{rq}
On ne peut pas supprimer l'hypothèse $\g(f,-1)^{\pp(e,0)}=\{0\}$. En effet, l'orbite $\Od_{1}$ de FII vérifie $\pp(e,0)=\{0\}$ et $\pp(e,1)\neq\{0\}$. Elle est cependant $\pp$-distinguée, donc $\pp$-self-large.\\
\end{Rq}

\begin{cor}
Les orbites $\Od_{50}$ de EV; $\Od_{85}$, $\Od_{88}$ de EVIII et $\Od_{16}$, $\Od_{17}$ de EI ne sont pas $\pp$-self-large.
\end{cor}
\begin{proof}
On a vu dans la section 6 que ces orbites nilpotentes sont presque $\pp$-distinguées et qu'elles ne sont pas paires. 
Grâce à \cite{Do1, Do2}, on peut calculer $\pp(e,0)$ et on trouve que $\pp(e,0)=\g(e,0)$ dans ces cas précis. On peut maintenant appliquer l'argument suivant de \cite{Pa4}.
Soit $\lf=\g^{\g(e,0)}$ et $\sfr=[\lf,\lf]$, de sorte que $\lf=\sfr\oplus\g(e,0)$. Alors $e$ est distingué dans $\sfr$ et la graduation induite par $h\in\sfr$ vérifie
$\{0\}=\sfr(-1,h)=\lf(-1,h)=\g(-1,h)^{\g(e,0)}=\g(-1,h)^{\pp(e,0)}\supseteq\g(f,-1)^{\pp(e,0)}$. 
En combinant ceci avec les tables de \cite{JN}, on montre que les hypothèses de la proposition \ref{sl} sont satisfaites. Les orbites mentionnées ne sont donc pas $\pp$-self-large.  
\end{proof}
En suivant la même idée, il est possible de donner une preuve alternative de la description de la section 5.2 des orbites presque $\pp$-distinguées de AI qui ne sont pas $\pp$-self-large.
Malheureusement, la proposition \ref{sl} ne permet pas l'étude des orbites presque $\pp$-distinguées non $\pp$-self-large de AII ni de l'orbite $\Od_{1}$ de EIV. 
Pour ces orbites, on se réfère aux sections 5.3 et 6.3. 
Finalement, les dernières orbites à traiter étant paires, on peut lister l'ensemble des orbites $\pp$-self-large dans les différents cas simples à l'aide du lemme \ref{norbite}.
Le cas général s'en déduit facilement étant donné que les algèbres de Lie symétriques sont produit direct d'algèbres de Lie symétriques simples.
\begin{center}
\begin{tabular}{|p{2cm}| p{10cm}|}
\hline
Cas &  Orbites $\pp$-self-large.\\
\hline
AI & L'orbite ($\pp$-distinguée) régulière et les orbites dont le diagramme associé est constituée de lignes de longueurs différant d'au moins 2.\\
\hline
AII & L'orbite ($\pp$-distinguée) régulière et les orbites dont le diagramme associé est constituée de paires de lignes de longueurs différant d'au moins 2.\\ 
\hline
AIII & Les orbites $\pp$-distinguées (\emph{i.e.} qui ont un $ab$-diagramme dont les lignes de même longueur débutent par la même lettre), ce sont les seules orbites presque $\pp$-distinguées.\\
\hline
BDI, CI & Les orbites presque $\pp$-distinguées (cf. section 3.3).\\
\hline
CII, DIII & Les orbites $\pp$-distinguées (qui sont les seules orbites presque $\pp$-distinguées, cf. section 3.3).\\
\hline
EIII, EVI, EVII, EIX, FI, FII, GI &  Les orbites $\pp$-distinguées (qui sont les seules orbites presque $\pp$-distinguées, cf. section 6.1).\\
\hline
EII & Les orbites presque $\pp$-distinguées. En particulier l'orbite $\Od_{22}$ non $\pp$-distinguée.\\
\hline
EIV & L'orbite régulière ($\pp$-distinguée).\\
\hline
EI & Les orbites $\pp$-distinguées et les orbites $\Od_{12}, \Od_{21}, \Od_{23}$ (cf. section 6 et lemme \ref{norbite}).\\
\hline
EV & Les orbites $\pp$-distinguées et l'orbite $\Od_{81}$ (cf. section 6 et lemme \ref{norbite}).\\
\hline
EVIII & Les orbites $\pp$-distinguées et les orbites $\Od_{81}, \Od_{95}$ (cf. section 6 et lemme \ref{norbite}).\\
\hline
Algèbres de Lie & Les orbites dont les éléments $e$ vérifient $\pp(e,0)$ est un tore et $\pp(e,1)=\{0\}$ (cf. \cite[Théoreme 2.1]{Pa4}).\\
\hline
\end{tabular}
\end{center}
\medskip

Les calculs permettent de montrer le fait suivant. Les éléments $\pp$-self-large d'une algèbres de Lie symétrique simple sont exactement les éléments $e$ vérifiant l'une des conditions suivantes
\begin{itemize}
\item $\pp(e,0)=\{0\}$ \emph{i.e.} $e$ est $\pp$-distingué;
\item $\pp(e,0)$ est un tore et $\pp(e,1)=\{0\}$. 
\end{itemize}
Cependant, ceci est uniquement valable dans le cas simple. En effet, la remarque \ref{rq} implique que, dans FII$\times$EI, 
l'orbite $\Od_{1}\times\Od_{21}$ est $\pp$-self-large mais vérifie $\pp(e,0)=T_{1}\neq \{0\}$ et $\pp(e,1)\neq \{0\}$.

\end{document}